%% file: main.tex
\documentclass[twocolumn]{autart}    %

\usepackage{algorithmic}
\usepackage{amsfonts}
\usepackage{amsmath}
\usepackage{amssymb}
\usepackage{anyfontsize}  
\usepackage{bbm, bm} 
\usepackage{cite}
\usepackage[dvipsnames]{color, xcolor} %
\usepackage{epsfig} %
\usepackage[T1]{fontenc}
\usepackage{graphicx}
\usepackage{mathrsfs}
\usepackage{pgfplots}
\usepgfplotslibrary{ternary, units}
\usepackage{subcaption}
\usepackage{svg}
\usepackage{textcomp}
\usepackage{tikz}
\usepackage{tikz-3dplot}
\usetikzlibrary{intersections,arrows,arrows.meta,shapes, positioning, shadows, shapes.geometric}
\usetikzlibrary{decorations.pathmorphing, pgfplots.ternary, pgfplots.units}
\usepackage{ulem}
\usepackage{url}
\usepackage[hidelinks]{hyperref}
\hypersetup{
   colorlinks   = true, %
   urlcolor     = blue, %
   linkcolor    = blue, %
   citecolor   = blue %
}
\edef\endfrontmatter{%
  \unexpanded\expandafter{\endfrontmatter}%
  \noexpand\endNoHyper %
}

\input{commands.sty}

\pgfplotsset{compat=1.7}
\begin{document}

\begin{frontmatter}

\title{Passivity Tools for Hybrid Learning Rules in Large Populations} %

\thanks[footnoteinfo]{
This work was supported in part by 
AFOSR Grant FA9550-23-1-0467 and 
NSF Grants 2139713, 2139781, 2139982,  
2135561, 
and 1846524.}%
\author[Maryland]{Jair Cert\'{o}rio}\ead{certorio@umd.edu}, { \author[USC]{Kevin Chang}\ead{kcchang@usc.edu},  
\author[Maryland]{Nuno C. Martins}\ead{nmartins@umd.edu}}, \author[USC]{Pierluigi Nuzzo}\ead{nuzzo@usc.edu},  \author[UC_Irvine]{Yasser Shoukry}\ead{yshoukry@uci.edu}  %

\address[Maryland]{Department of ECE and ISR, University of Maryland, College Park, MD, 20742, USA.}
\address[USC]{Department of Electrical and Computer Engineering, University of Southern California, Los Angeles, CA, 90089, USA.}
\address[UC_Irvine]{Department of Electrical Engineering and Computer Science, University of California, Irvine, CA, 92697, USA.}


\begin{abstract}                          %
 Recent work has pioneered the use of system-theoretic passivity to study equilibrium stability for the dynamics of noncooperative strategic interactions in large populations of learning agents.  In this and related works, the stability analysis leverages knowledge that certain ``canonical'' classes of learning rules used to model the agents' strategic behaviors satisfy a passivity condition known as $\delta$-passivity. In this paper, we consider that agents exhibit learning behaviors that do not align with a canonical class. Specifically, we focus on characterizing $\delta$-passivity for hybrid learning rules that combine elements from canonical classes. Our analysis also introduces and uses a more general version of $\delta$-passivity, which, for the first time, can handle discontinuous learning rules, including those showing best-response behaviors. We state and prove theorems establishing $\delta$-passivity for two broad convex cones of hybrid learning rules. These cones can merge into a larger one preserving $\delta$-passivity in scenarios limited to two strategies. In our proofs, we establish intermediate facts that are significant on their own and could potentially be used to further generalize our work. We illustrate the applicability of our results through numerical examples.
\end{abstract}

\end{frontmatter}

\section{Introduction}

 Developments in population games and evolutionary dynamics~\cite{Sandholm2010Population-Game,Sandholm2015Handbook-of-gam} have contributed to systematic methods to model and analyze the dynamics of strategic noncooperative interactions among large populations of learning agents. Central to this approach is the use of learning rules\footnote{Often called revision protocols in the population game literature.} that model how agents update their strategies over time based on the payoffs of those strategies, which are in turn determined by a payoff mechanism. These learning rules can be explicitly programmed in artificial agents or represent innate preferences or bounded rationality in humans or other natural agents. This framework is well-suited to distributed optimization~\cite{Barreiro-Gomez2016Constrained} and engineering systems~\cite{Tembine2010Evolutionary-ga,Obando2014Building-temper}, where the payoff mechanisms might abstractly represent the specifics of the agents' interaction environments, such as in congestion games, or are purposefully designed and implemented by a coordinator to steer the population toward desired strategic outcomes~\cite{Quijano2017The-role-of-pop,Martins2023Epidemic-popula}.

\subsection{Studying Passivity for Hybrid Learning Rules}
System-theoretic passivity tools have been employed to study how populations, adhering to certain learning rules, achieve and maintain Nash equilibria---a process often referred to as Nash equilibrium seeking. Pioneering work in~\cite{Fox2013Population-Game} demonstrates that Nash equilibrium seeking and convergence are achieved under contractive payoff mechanisms if the learning rules are $\delta$-passive, a concept inspired by classical notions of system-theoretic passivity that has been generalized in~\cite{Park2019From-Population,Arcak2020Dissipativity-T}. The convergence results in these articles allow for dynamic payoff mechanisms and are not contingent upon the specific learning rule employed, as long as it satisfies $\delta$-passivity, often confirmed through structural analysis.

Recent research has shown that $\delta$-passivity can also be used to systematically design dynamic payoff mechanisms to guarantee, via Lyapunov analysis, global asymptotic stability of optimal strategic profiles in the presence of coupled dynamics and constraints.

Therefore, it is understandable that ongoing work focuses on broadening well-established classes of $\delta$-passive learning rules. The broader these classes are, the more universally applicable the results are for $\delta$-passive learning rules, thus decreasing the need for detailed knowledge of specific rules. In real-world scenarios where agents exhibit behaviors that do not align perfectly with any existing class, it is essential to go further and allow for learning rules outside known $\delta$-passive classes. Consequently, we focus on characterizing $\delta$-passivity for the so-called hybrid learning rules that do not belong to an existing class but combine elements from well-established classes.

There is an alternative method employing equilibrium-independent passivity~\cite{Hines2011Equilibrium-ind} to characterize equilibrium stability in population games for replicator dynamics~\cite{Mabrok2021Passivity-Analy}. Similar equilibrium-independent passivity techniques have been applied to establish convergence properties for a modification of replicator dynamics that incorporates exponentially-discounted learning~\cite{Gao2020On-Passivity-Re}. Unfortunately, these methods do not apply to our framework because equilibrium-independent passivity has not been established for the evolutionary dynamics associated with the canonical or hybrid classes of learning rules considered here. Conversely, \cite[Proposition~III.5]{Park2018Passivity-and-e} demonstrates that replicator dynamics is not $\delta$-passive, making our approach unsuitable for replicator dynamics.  

\subsection{Contributions and Structure of the Paper}

This paper is the first to characterize $\delta$-passivity for hybrid learning rules. Although we do not fully characterize $\delta$-passivity for all possible hybrid learning rules, we make substantial progress. We establish $\delta$-passivity of two broad convex cones of hybrid learning rules. For the case of two strategies, we demonstrate that these convex cones can be further combined into a larger cone that also maintains $\delta$-passivity. Our research further introduces a generalization of $\delta$-passivity that accommodates the inclusion of the so-called best-response learning rule into the convex cones of hybrid learning rules, enabling the handling of discontinuous dynamics previously incompatible with traditional $\delta$-passivity approaches.

After the Introduction, \S\ref{section2} describes our framework, in which the payoff mechanism interacts in feedback with an evolutionary dynamics model derived to capture the effect of the learning rule. Moreover, we incorporate Filippov solution concepts to manage discontinuous dynamics. In \S\ref{sec:DeltaPassivityAndNES} we generalize existing $\delta$-passivity concepts, list the canonical classes of learning rules we intend to merge into hybrid rules, and explore the importance of $\delta$-passivity. 
In \S\ref{sec:HybridLearnRuleMainProb} we define these hybrid learning rules and formulates our main problem. Our main results are presented in \S\ref{sec:MainResults}, where we provide conditions to maintain $\delta$-passivity and illustrate examples of hybrid rules that adhere to these conditions. 
In \S\ref{sec:AuxResults} we provide additional propositions to prove our main results while \S\ref{sec:sims} illustrates them via simulations on the case study of a weighted congestion game. Finally, \S\ref{sec:conclusion} summarizes our results and future work. 

\section{Framework and Problem Formulation}
\label{section2}

Our framework models the noncooperative strategic interactions of a large number of agents. For simplicity, we assume that all agents belong to one population presented with a finite set of available strategies $\{1,\ldots,n\}$.
Each agent follows one strategy at a time, which the agent can revise repeatedly. A payoff vector $p\st$ in $\mathbb{R}^n$, whose $i$-th entry $p_i\st$ quantifies the payoff (net reward) of the strategy $i$ at time $t$, influences the revision process, as agents typically seek strategies with higher payoffs. Agents are non-descript, and the so-called population state vector $x\st$ in ${\mathbb{X} := \{ x \in \mathbb{R}^n \ | \ x_i \geq 0, \ \sum_{i=1}^n x_i = 1 \}}$ approximates the population's aggregate strategic choices in the large population limit~\cite{Sandholm2003Evolution-and-e,Kara_2023aa}. That is, $x_i\st$ approximates the proportion of agents in the population following strategy $i$ at time $t$. 

\paragraph*{Notation.} Throughout this paper, $\mathbb{R}^n$ and $\mathbb{R}_{\geq0}^n$ denote the set of $n$-dimensional real column vectors and nonegative column vectors, respectively. Let $\mathbf{1}$ denote the vector of ones and $I$ denote the identity matrix, whose dimension will be clear from the context. Let $v'$ denote the transpose of $v$, and $A \succeq B$ denote that $A - B$ is positive semidefinite. 
For $x\in \mathbb{R}^n$, $\|x\|$ denotes an arbitrary norm, and is used for statements that hold independently of the norm used, while $\|x\|_p$ denotes a particular $p$-norm. Further notation will be introduced as needed.

\subsection{Learning Rules}

A piecewise continuous map $\mathcal{T}: \mathbb{X} \times \mathbb{R}^{n} \rightarrow \mathbb{R}_{\geq 0}^{n \times n}$,  %
referred to as {\it learning rule}\footnote{Learning rules, as we define here, are also denoted in the literature as revision protocols.}, models the agents' strategy revision preferences. Specifically, $\mathcal{T}_{ij}(x,p)$ is the rate at which agents following strategy $i$ switch to strategy $j$ as a function of $x$ and $p$. In \cite[Part~II]{Sandholm2010Population-Game} and \cite[\S 13.3-13.5]{Sandholm2015Handbook-of-gam} there is a comprehensive discussion on the types of learning rules and the classes of bounded rationality behaviors that they model.

\begin{remark}\label{rem:ContPropLearn} {\bf (Continuity properties for $\mathcal{T}$)} As will be clear from our analysis, any discontinuity in $\mathcal{T}$ will result from incorporating the so-called {\it best response} learning rule, which we will introduce in~\S\ref{sec:KeyLearningRules}.  We will consider learning rules satisfying the following properties:
\begin{enumerate}
    \item  $\mathcal{T}$ is locally essentially bounded.
    \item For any given $p$ in $\mathbb{R}^n$, $\mathcal{T}(x,p)$ is Lipschitz continuous with respect to $x$.
    \item For any given $x$ in $\mathbb{X}$, $\mathcal{T}(x,p)$ is piecewise continuous with respect to $p$. All points $p$ at which $\mathcal{T}(x,p)$ is discontinuous for some $x \in \mathbb{X}$ are contained in the zero-measure set 
    $$ \mathbb{D}:=\{ p \in \mathbb{R}^n \ | \ p_i=p_j \text{ for some $i \neq j$} \}.$$
\end{enumerate}
\end{remark}

The revision process causes $x$ to vary over time. Following the standard approach in~\cite[\S4.1.2]{Sandholm2010Population-Game}, the {\it evolutionary dynamics model} below governs the dynamics of~$x\st$:
\begin{subequations}
\begin{equation}\tag{EDMa} \label{eq:EDM-DEF} \dot {x}\st \underset{p(t) \notin \mathbb{D}}{=} \mathcal{V} ( x\st,p\st ), \quad t\geq 0, 
\end{equation}
where the $i$-th component of $\mathcal{V}$ accounts for the net flow of agents switching to strategy $i$ according to 
\begin{equation} \tag{EDMb} \label{eq:EDM-DEF-b}
\mathcal{V}_i ( x,p ) := \sum_{j=1}^n \underbrace{x_j \mathcal{T}_{ji}(x,p)}_{flow \ j\rightarrow i} - \underbrace{x_i \mathcal{T}_{ij}(x,p)}_{flow \ i\rightarrow j}.
\end{equation}
\end{subequations}  

In order to guarantee the existence of solutions later, we do not require the equality in~(\ref{eq:EDM-DEF}) to hold when $p\st$ is in the discontinuity set $\mathbb{D}$. To be specific, we will adopt a standard notion of a solution based on Filippov's differential inclusion method explained in~\cite[p.48]{Cortes_2008aa}. We refer to the evolutionary dynamic model and its governing equations~\eqref{eq:EDM-DEF}-\eqref{eq:EDM-DEF-c} by (EDM). Equation~\eqref{eq:EDM-DEF-c} will be introduced in the following section, along with the solution concept.

\subsection{Payoff Mechanism and Filippov Solutions}
\label{sec:PayoffMechAndFilippovSoln}

\begin{figure}
\begin{center}
\begin{tikzpicture}[scale=1.0,
  transform shape,
  node distance=1.5cm,
  block/.style={rectangle, draw, rounded corners, fill=teal!10, minimum width=3cm, minimum height=1cm, align=center},
  arrow/.style={->, >=stealth, thick}
  ]

  \node[block, fill=yellow!10] (strat) {\Large $\dot{x}=\mathcal{V}(x,p)$};
  \node[block, below=of strat] (payoff) {\Large $\mathcal{F}$ or PDM};
  \node[yshift=-23] at (payoff) {\it payoff mechanism};

  \draw[arrow] (strat) -- (payoff) node[midway, right] {\Large $x$};
  \draw[arrow] (payoff.east) -- ++(2,0) |- (strat.east) node[midway, above] {\Large $p$};

\end{tikzpicture}
\end{center}
\caption{Interconnection of (EDM) and payoff mechanism.}
\label{fig:closedloop}
\end{figure}
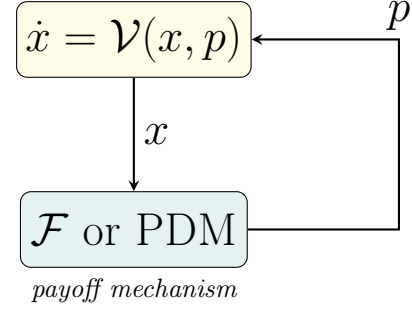

Before we describe our solution concept, we have to start by stipulating that a payoff mechanism is a causal map that takes ${x:[0,\infty) \rightarrow \mathbb{X}}$ as input and produces ${p:[0,\infty) \rightarrow \mathbb{R}^n}$ as its output. Hence, the payoff mechanism and (EDM) are interconnected in feedback as in Fig.~\ref{fig:closedloop} and our notion of solution will be with respect to the initial value problem associated with the state of the entire feedback system, which will be $x:[0,\infty) \rightarrow \mathbb{X}$ if the payoff mechanism is memoryless and ${(x,q): [0,\infty) \rightarrow \mathbb{X} \times \mathbb{R}^m}$ if the payoff mechanism has memory with a finite-dimensional state $q:[0,\infty) \rightarrow \mathbb{R}^m$.

\begin{assumption}{\bf (Valid $p:[0,\infty) \rightarrow \mathbb{R}^n$)} We restrict our analysis to payoff mechanisms producing ${p:[0,T] \rightarrow \mathbb{R}^n}$ absolutely continuous for any $T \geq 0$. We refer to any $p:[0,\infty) \rightarrow \mathbb{R}^n$ satisfying this condition as \underline{valid}.
\end{assumption}

\begin{definition}{\bf (Filippov solutions\cite{Cortes_2008aa})} \label{def:Filippov} We start by adapting~\cite[equation (19)]{Cortes_2008aa} to our framework by defining the set-valued map $\mathfrak{V}:\mathbb{X}\times \mathbb{R}^n \rightarrow 2^{\mathbb{R}^n}$ specified as
\begin{equation}
\mathfrak{V}(x,p) := 
\cap_{ \delta > 0} \bar{co} \big \{ \mathcal{V}(x,p + \tilde{p}) \ | \ \| \tilde{p}\| \leq \delta  \big \},
\end{equation} where $\bar{co}$ produces the convex hull of a set.
In our analysis, ${x:[0,\infty) \rightarrow \mathbb{X}}$ is a solution for~\eqref{eq:EDM-DEF} in the sense of Carath\'{e}odory, which is required to satisfy the differential inclusion\footnote{Notice that $\mathfrak{V}(x,p) = \{ \mathcal{V}(x,p) \}$ when $p \notin \mathbb{D}$.}
\begin{equation}
\tag{EDMc} \label{eq:EDM-DEF-c}
\dot{x}\st \in \mathfrak{V}(x\st,p\st), \quad t \geq 0
\end{equation} for almost every $t$ and ${x:[0,T] \rightarrow \mathbb{X}}$ is absolutely continuous for every $T >0$. We adopt the widely-used convention of referring to any ${x:[0,\infty) \rightarrow \mathbb{X}}$ obtained in this manner as a \underline{Filippov solution} for~\eqref{eq:EDM-DEF}.
\end{definition}

\begin{definition}{\bf (Game)} A payoff mechanism is referred to as a {\it game} when it is a memoryless map ${\mathcal{F}:x\st \mapsto p\st}$, where $\mathcal{F}:\mathbb{X} \rightarrow \mathbb{R}^n$ is Lipschitz continuous. 
\end{definition}

\begin{remark} We will argue that for payoff mechanisms of interest, the Filippov solutions for the system shown in Fig.~\ref{fig:closedloop} ensure that the function \( p:[0,\infty) \to \mathbb{R}^n \) is valid. Since we know from Definition~\ref{def:Filippov} that \( {x:[0,T] \to \mathbb{X}} \) is absolutely continuous for any \( T>0 \), confirming that \( {p:[0,\infty) \to \mathbb{R}^n} \) is valid means that \underline{both \( x \) and \( p \) } functions are \underline{almost everywhere} \underline{differentiable} over time.
\end{remark}

If the payoff mechanism is a game $\mathcal{F}$, then the closed-loop system in Fig.~\ref{fig:closedloop} is ${\dot {x}\st = \mathcal{V} ( x\st,\mathcal{F}(x\st) )}$, for all $t$ for which $\mathcal{F}(x\st) \notin \mathbb{D}$. Invoking Remark~\ref{rem:ContPropLearn}, \cite[Proposition~3]{Cortes_2008aa} applied to this system guarantees that there is a Filippov solution $x: [ 0,\infty) \rightarrow \mathbb{X}$ for any $x\so$ in $\mathbb{X}$. Since $\mathcal{F}$ is Lipschitz continuous, the absolute continuity of $x:[0,T] \rightarrow \mathbb{X}$, for any $T>0$, guarantees the validity of $p:[0,\infty) \rightarrow \mathbb{R}^n$.

More broadly, the payoff mechanism might also be a finite-dimensional dynamical system, as described in~\cite{Fox2013Population-Game} and explored in detail in~\cite{Park2019From-Population}, where it is called a payoff dynamical model (PDM). If the payoff mechanism is a PDM and its state at time $t$ is $q\st$, then the closed-loop system illustrated in Fig.~\ref{fig:closedloop}, with state $(x\st,q\st)$ at time $t$, is piecewise continuous and meets the conditions of \cite[Proposition~3]{Cortes_2008aa}. Thus, Filippov solutions exist for any initial state $(x\so,q\so)$ in $\mathbb{X} \times \mathbb{R}^n$. From \cite[Definition~7]{Park2019From-Population}, it follows by an argument  similar to the one above
that the function $p:[0,\infty) \rightarrow \mathbb{R}^n$ is valid.

\section{$\delta$-Passivity and Nash Equilibria Seeking}
\label{sec:DeltaPassivityAndNES}

A Lyapunov method was used in~\cite{Hofbauer2009Stable-games-an} to prove global asymptotic convergence of $x\st$ to the set of Nash equilibria of the so-called contractive games. This convergence is often desirable and is commonly referred to as {\it Nash equilibria seeking}. The notion of $\delta$-passivity was proposed in~\cite{Fox2013Population-Game} as a means of generalizing the Lyapunov method in~\cite{Hofbauer2009Stable-games-an} so as to allow for more general payoff mechanisms, including PDMs. Subsequent work in~\cite{Park2019From-Population,Kara2021Pairwise-Compar} significantly expanded the application of $\delta$-passivity to important classes of learning rules. 

Previous definitions of $\delta$-passivity are incompatible with discontinuous learning rules, such as the best response learning rule we will define in \S\ref{subsec:DeltaPassForDesign}. Consequently, we start with the following slightly different definition of $\delta$-passivity, which is general enough to expand its use to the classes of discontinuous hybrid learning rules we will investigate later in \S\ref{sec:HybridLearnRuleMainProb}.

\begin{definition} \label{def:deltapassivity} \textbf{($\delta$-Passivity)} Given a learning rule $\mathcal{T}$, the corresponding (EDM) is said to be $\delta$-passive when there is a locally Lipschitz function $\mathcal{S}:\mathbb{X} \times \mathbb{R}^n \rightarrow [0,\infty)$ and a continuous function $\mathcal{P}:\mathbb{X} \times \mathbb{R}^n \rightarrow [0,\infty)$ such that for any given valid $p:[0,\infty) \rightarrow \mathbb{R}^n$ and any ${x:[0,\infty) \rightarrow \mathbb{X}}$ obtained from a Filippov solution, the following inequality is satisfied for almost all $t \geq 0$
\begin{subequations}
\label{eq:DeltaPassDef}
\begin{equation}
\label{eq:DeltaPassDef(a)}
\dot{s}\st \leq  \ \dot{x}(t)'\dot{p}\st - \wp\st,
\end{equation} where $s\st : = \ \mathcal{S}(x\st,p\st)$ and $\wp\st:=\mathcal{P}(x\st,p\st)$. In addition, the following equivalence must hold for all $\ x \in \mathbb{X}$ and  $p \in \mathbb{R}^n$
\begin{equation}
\label{eq:DeltaPassDef(b)}
\mathcal{P}(x,p) = 0 \Leftrightarrow \mathcal{S}(x,p)=0 \Leftrightarrow x \in \mathfrak{B}(p),
\end{equation} where $\mathfrak{B}(p)$ is the best response set-valued map
\end{subequations}
$$\mathfrak{B}(p) := \Big \{z \in \mathbb{X} \ \Big | \ p'z = \max_{1\leq i \leq n} p_i \Big \}.$$
\end{definition}

Notice that since $\mathcal{S}$ is locally Lispchitz, absolute continuity of both $x:[0,T] \rightarrow \mathbb{X}$ and $p:[0,T] \rightarrow \mathbb{R}^n$, for any $T \geq 0$, guarantees that $s:[0,\infty) \rightarrow [0,\infty)$ is differentiable almost everywhere. Consequently, $x\st$, $p\st$ and $s\st$ are \underline{simultaneously differentiable} for almost all $t \geq 0$ and is at these values of $t$ that~(\ref{eq:DeltaPassDef(a)}) must hold to ensure $\delta$-passivity.

\subsection{Comparing Definition~\ref{def:deltapassivity} with Standard $\delta$-Passivity}

There are two main differences between Definition~\ref{def:deltapassivity} and the standard definitions of $\delta$-passivity in~\cite{Fox2013Population-Game,Park2019From-Population,Arcak2020Dissipativity-T}. The first difference is allowing $\mathcal{S}$ to be locally Lipschitz rather than continuously differentiable as in~\cite{Fox2013Population-Game}. The second difference is replacing the original equivalence requirement $\mathcal{S}(x,p)=0 \Leftrightarrow \mathcal{V}(x,p) = 0$ in~\cite{Arcak2020Dissipativity-T} with $\mathcal{S}(x,p)=0 \Leftrightarrow x \in \mathfrak{B}(p)$ as we do in~(\ref{eq:DeltaPassDef(b)}). Both of these equivalence requirements are interchangeable for the large classes of learning rules considered in~\cite{Fox2013Population-Game,Park2019From-Population,Arcak2020Dissipativity-T,Kara2021Pairwise-Compar}, which are referred to as {\it Nash stationary}~\cite{Sandholm2015Handbook-of-gam}. However, these equivalence requirements are not interchangeable for the more general classes of learning rules we now consider, including the best response. By imposing $\mathcal{S}(x,p)=0 \Leftrightarrow x \in \mathfrak{B}(p)$ in~(\ref{eq:DeltaPassDef(b)}), as we illustrate in~\S\ref{subsubsec:deltapassivitycontractivegames}-\S\ref{subsec:DeltaPassForDesign}, one can establish Nash equilibria seeking or other analogous properties indirectly by showing that $s\st$ vanishes as $t$ tends to infinity.

\subsection{Canonical Classes of $\delta$-Passive Learning Rules}
\label{sec:KeyLearningRules}

We proceed to define three well-known and very general learning protocol classes that are  $\delta$-passive.

\begin{definition} {\bf (Best response learning rule)} \label{def:BRLR}
The so-called best response protocol is defined as
\begin{align}
\mathcal{T}^{\text{\tiny BR}}_{ij}(x,p)&=\mathcal{Y}_j(p), \label{eqdef:BR} \\ \nonumber  \mathcal{Y}_j(p) & := \lim_{\beta \rightarrow 0^+} \frac{e^{p_j / \beta}}{\sum_{\ell=1}^n e^{p_\ell / \beta}},
\end{align} for all $x$ in $\mathbb{X}$ and $p$ in $\mathbb{R}^n$. 
\end{definition}

The parameter $\beta$ can be viewed as the noise parameter of the so-called logit learning rule\footnote{See \cite[\S2]{Hofbauer2002On-the-global-c} for a discussion.}, in which case its noise-free limit is $\mathcal{T}^{\text{\tiny BR}}$.

Here, $\mathcal{Y}(p)$ is the limit as $\beta$ tends to zero of the gradient of the softmax function $\beta \ln(e^{p_1 / \beta} + \cdots + e^{p_n / \beta})$. Its entries can be written explicitly as $\mathcal{Y}_i(p) = 0$ if there is a strategy $\ell$ for which $p_i < p_\ell$, and, otherwise, $\mathcal{Y}_i(p) = m(p)^{-1}$, where $m(p)$ is the number of strategies with maximal payoff, that is, $m(p)$ is the number of elements of 
$$\mathbb{M}(p):={ \Big \{ 1 \leq i \leq n  \ \Big | \ p_i = \max_{1 \leq \ell \leq n} p_\ell \Big  \}}.$$

Hence, if an agent following the best response learning rule switches strategies, it does so to a uniformly randomly selected strategy in $\mathbb{M}(p)$. 

\begin{remark}[\bf Justification for Definition~\ref{def:BRLR}] \label{rem:ComparisonStandardDefBR} 

In its standard definition, such as in~\cite[\S13.5.2.2]{Sandholm2015Handbook-of-gam}, the best response learning rule is a set-valued map taking values in $\mathfrak{B}(p)$, leading to the differential inclusion~\cite[(13.17)]{Sandholm2015Handbook-of-gam} instead of (EDM). The standard definition and ours in~(\ref{eqdef:BR}) are identical for almost all $p$ because $\mathfrak{B}(p) = \{ \mathcal{Y}(p) \}$ for all $p$ in $\mathbb{R}^n-\mathbb{D}$. 
Furthermore, for any given payoff mechanism that is a game $\mathcal{F}$, one can show that the differential inclusion~\cite[(13.17)]{Sandholm2015Handbook-of-gam} is identical to that constructed as in \cite[(19)]{Cortes_2008aa} for $\dot{x} = \mathcal{V}^{\text{\tiny BR}}(x,\mathcal{F}(x))$, where $\mathcal{V}^{\text{\tiny BR}}$ results from~(\ref{eqdef:BR}) substituted into (EDM). Hence, by construction, the set of Filippov solutions for $\dot{x} = \mathcal{V}^{\text{\tiny BR}}(x,\mathcal{F}(x))$ is identical to the set of Carath\'{e}odory solutions of~\cite[(13.17)]{Sandholm2015Handbook-of-gam}. More generally, for any given payoff mechanism considered in~\S\ref{sec:PayoffMechAndFilippovSoln}, the set of Filippov solutions for the closed-loop system in Fig.~\ref{fig:closedloop} assuming~(\ref{eqdef:BR}) is identical to the set of Carath\'{e}odory solutions of the differential inclusions we would have obtained if we had used the standard set-valued definition for the best response learning rule. Our option to use~(\ref{eqdef:BR}) along with Filippov solutions is, then, interchangeable with the standard approach but has the benefit of simplifying our notation by avoiding the need to consider differential inclusions when analyzing the hybrid learning rules that will be defined later on in~\S\ref{sec:HybridLearnRuleMainProb}.
\end{remark}

We can now invoke Remark~\ref{rem:ComparisonStandardDefBR} and conclude, by replacing $\mathcal{F}(x\st )$ with a valid $p\st$, that the analysis in~\cite[Appendix~A.3.]{Hofbauer2009Stable-games-an} leads directly to the following proposition.

\begin{proposition}\label{prop:BRDeltaPassive} The best response learning rule is $\delta$-passive with storage function
$$ \mathcal{S}^{\text{\tiny BR}}(x,p)= \mathcal{P}^{\text{\tiny BR}}(x,p): = p'\mathcal{Y}(p) - p'x, \quad x \in \mathbb{X}, \ p \in \mathbb{R}^n.$$
\end{proposition}

Proposition~\ref{prop:BRDeltaPassive} ascertains explicitly, for the first time, that the best-response learning rule is $\delta$-passive. 
{
\begin{definition}{\bf ($\mathbb{T}^\text{\tiny IPC}$)} \label{def:IPC} Any protocol is of the {\it impartial pairwise comparison} (IPC) type~\cite{Sandholm2010Pairwise-compar} if for each $j \in [n]$ there is a map ${\psi}_j:\mathbb{R}_{\geq0}\rightarrow \mathbb{R}_{\geq0}$, such that ${\psi}_j(\nu)>0$ for $\nu>0$ and ${\psi}_j(0)=0$, and $\cT$ can be recast as:
\begin{equation}
\label{eq:defIPC}
\mathcal{T}_{ij}(x,p) \underset{\text{\tiny IPC}}{=} {\psi}_j([\tilde{p}_{ij}]_+),
\end{equation} where $\tilde{p}_{ij}:=p_j-p_i$. We use $\mathbb{T}^\text{\tiny IPC}$ to denote the set of all IPC learning rules.
\end{definition}

An important subset of the IPC rules is the one for which the agents' switching rate from strategy $i$ to $j$ does not decrease as the payoff difference $\tilde{p}_{ij}$ increases. Such set contains the most sensible rules within $\mathbb{T}^\text{\tiny IPC}$.

\begin{definition}{\bf ($\mathbb{T}^\text{\tiny ND}$)} \label{def:ND} We denote by $\mathbb{T}^\text{\tiny ND}$ the set of nondecreasing IPC (ND-IPC) learning rules, that are defined by a $\psi$ such that $\psi_j$ is nondecreasing for all $j$.
\end{definition}
}

\begin{example} \label{ex:Smith} The well-known Smith's learning rule~\cite{Smith1984The-stability-o} is of {the ND-IPC class} and can be specified by selecting ${\psi_j^{\text{\tiny Smith}}([\tilde{p}_{ij}]_+) := \lambda [\tilde{p}_{ij}]_+}$, ${\lambda>0}$. 
\end{example}

\begin{proposition} \label{prop:IPCDeltaPassive}
Any IPC learning rule is $\delta$-passive, where for all $x \in \mathbb{X}$ and $ p \in \mathbb{R}^n$, we can use
\begin{align} \mathcal{S}^{\text{\tiny IPC}}(x,p) : = &\sum_{i=1}^n\sum_{j=1}^n x_i \int_0^{\tilde{p}_{ij}}\psi_j(\tau)d \tau, \\ \mathcal{P}^{\text{\tiny IPC}}(x,p) = &-\sum_{i=1}^n\sum_{j=1}^n \mathcal{V}_i(x,p) \int_0^{\tilde{p}_{ij}}\psi_j(\tau)d \tau.
\end{align}
\end{proposition}

{
\begin{definition}{\bf ($\mathbb{T}^\text{\tiny SEPT}$)} \label{def:SEPT} Any protocol is of the {\it separable excess payoff target} (SEPT) type~\cite{Sandholm2005Excess-payoff-d} if for each $j \in [n]$ there is a map ${\phi}_j:\mathbb{R}_{\geq0}\rightarrow \mathbb{R}_{\geq0}$, satisfying ${\phi}_j(\nu)>0$ for $\nu>0$ and ${\phi}_j(0)=0$, and $\cT$ can be recast as:
\begin{equation}
\label{eq:defSEPT}
\mathcal{T}_{ij}(x,p) \underset{\text{\tiny SEPT}}{=} {\phi}_j([\hat{p}_{j}]_+), \quad \hat{p}_{j}:= p_j-\sum_{i=1}^n x_i p_i,
\end{equation} where $[\hat{p}_{j}]_+ := \max(0,\hat{p}_{j})$ and $\hat{p}$ is the so-called excess payoff vector. Let $\mathbb{T}^\text{\tiny SEPT}$ denote the set of all IPC rules.
\end{definition}}

\begin{example} \label{ex:BNN} The classic Brown-von Neumann-Nash (BNN) learning rule~\cite{Brown1950Solutions-of-ga} is of the SEPT class and is specified by $\phi_j^{\text{\tiny BNN}}([\hat{p}_{j}]_+) := \lambda [\hat{p}_{j}]_+$, ${\lambda>0}$.
\end{example}

\begin{proposition} \label{prop:SEPTDeltaPassive}
Any SEPT learning rule is $\delta$-passive, where for all $x \in \mathbb{X}$ and $ p \in \mathbb{R}^n$, we can use
\begin{align} \mathcal{S}^{\text{\tiny SEPT}}(x,p) : = &\sum_{j=1}^n \int_0^{[\hat{p}_{j}]_+}\phi_j(\tau)d \tau, \\ \mathcal{P}^{\text{\tiny SEPT}}(x,p) = & \sum_{i=1}^n\sum_{j=1}^n \phi_i([\hat{p}_i]_+) p_j \mathcal{V}_j(x,p) .
\end{align}
\end{proposition}

\subsection{$\delta$-Passivity and Contractive Games}
\label{subsubsec:deltapassivitycontractivegames}

 To briefly illustrate how to use a storage function satisfying~(\ref{eq:DeltaPassDef}) to study global asymptotic stability, consider that $p\st = \mathcal{F} (x \st)$ for a continuously differentiable contractive game\footnote{Contractive games were introduced in~\cite{Hofbauer2009Stable-games-an}, where they were originally called stable games.} $\mathcal{F}$. The (EDM) in this case specializes to $\dot{x}\st=\mathcal{V}(x\st,\mathcal{F}(x\st))$ and the following inequality holds for any Filippov solution and for all $t \geq 0$ where $x\st$ is differentiable
\begin{equation} 
\label{eq:ContractiveGameDerStorage}
\dot{s}\st \leq - \wp\st.
\end{equation}

In order to arrive at (\ref{eq:ContractiveGameDerStorage}) we invoke Remark~\ref{rem:ComparisonStandardDefBR} and follow the same steps\footnote{ As we explain in Remark~\ref{rem:ComparisonStandardDefBR}, the Carath\'{e}odory solution set in~\cite[Appendix~A.3.]{Hofbauer2009Stable-games-an} is identical to the Filippov solution set considered here.} as in~\cite[Appendix~A.3.]{Hofbauer2009Stable-games-an}. The analysis in~\cite[Appendix~A.3.]{Hofbauer2009Stable-games-an} considers the best response learning rule and also holds trivially for any $\mathcal{T} \in \mathbb{T}^{\text{\tiny IPC}} \cup \mathbb{T}^{\text{\tiny SEPT}}$ or, more generally, any $\delta$-passive learning rule.%

Since $s$ is differentiable for almost all $t \geq 0$, (\ref{eq:ContractiveGameDerStorage}) implies that $\lim_{t \rightarrow \infty} s\st=\lim_{t \rightarrow \infty} \wp\st = 0$ and consequently Nash equilibria seeking succeeds. Specifically, from~(\ref{eq:DeltaPassDef(b)}) we conclude that $x\st$ will converge globally asymptotically to the set $\mathbb{NE}(\mathcal{F})$ of Nash equilibria of $\mathcal{F}$ defined as
$$ \mathbb{NE}(\mathcal{F}) : = \{ x \in \mathbb{X} \ | \ x \in \mathfrak{B}(\mathcal{F}(x)) \}, $$ and the function $\mathcal{S}(x,\mathcal{F}(x))$ is a Lyapunov function.

\subsection{$\delta$-Passivity and Payoff Dynamic Models}
\label{subsec:DeltaPassivityUsefulPDM}

The framework and results in~\cite{Fox2013Population-Game} show that $\delta$-passivity and the existence of a storage function can be crucial for stability analysis when $p\st$ is generated by a PDM. A characterization of PDM classes for which stability results hold for $\delta$-passive learning rules was initially introduced in~\cite{Fox2013Population-Game} and later generalized in~\cite{Park2019From-Population}. This characterization was further extended for the more general notion of $\delta$-dissipativity in~\cite{Arcak2020Dissipativity-T}.

\subsection{Using $\delta$-Passivity for Design}
\label{subsec:DeltaPassForDesign}

The method used in~\cite{Martins2023Epidemic-popula,Certorio2022Epidemic-Popula} exemplifies how $\delta$-passivity can be used for mechanism design with performance guarantees, even when there are higher-order dynamics coupled to (EDM). It highlights the practical benefit gained from the facts that $\delta$-passivity is {\bf (i)} equilibrium independent (as defined in~\cite{Hines2011Equilibrium-ind}) and {\bf(ii)} the associated storage function can be used to construct a Lyapunov function establishing global asymptotic stability properties and performance bounds. Importantly, equilibrium independence allows for the decoupled design of {\bf (i)} a desirable equilibrium and {\bf (ii)} the dynamics of the payoff mechanism that will render the equilibrium globally asymptotically stable.   

\section{Hybrid Learning Rules and Main Problem}
\label{sec:HybridLearnRuleMainProb}

The learning rule $\mathcal{T}^\text{\tiny BR}$ and the classes $\mathbb{T}^\text{\tiny IPC}$ and $\mathbb{T}^\text{\tiny SEPT}$ model different preferences in strategic behaviors of agents. Specifically, agents following $\mathcal{T}^\text{\tiny BR}$ are quite aggressive and only adopt strategies with the highest payoffs. In contrast, a rule within $\mathbb{T}^\text{\tiny SEPT}$ allows switching to strategies that simply exceed the average population payoff, while rules in $\mathbb{T}^\text{\tiny IPC}$ allow agents to choose strategies that improve their current payoffs. Therefore, these rules encompass behaviors ranging from pursuing the maximum payoff to enhancing current strategy payoffs.

These learning classes are also quite broad, as they are \emph{convex cones}; this implies that they encompass any positive linear combination of their elements. Nevertheless, it is unrealistic to expect that a human or other natural agent will adhere strictly to a rule from one of these classes. A more realistic scenario, as discussed by Sandholm \cite[p.~5]{Sandholm2010Pairwise-compar}, suggests that an agent's learning rule might be hybrid, incorporating elements from multiple classes as specified below.

 \begin{definition}{\bf (\bf Hybrid learning rule)} \label{def:HybLearnRule} We qualify a learning rule $\mathcal{T}$ as hybrid if it is expressible as $$\mathcal{T} = \alpha^{\text{\tiny BR}} \mathcal{T}^{\text{\tiny BR}} + \alpha^{\text{\tiny SEPT}}\mathcal{T}^{\text{\tiny SEPT}}+\alpha^{\text{\tiny IPC}}\mathcal{T}^{\text{\tiny IPC}},$$ where $\alpha^{\text{\tiny BR}}$, $\alpha^{\text{\tiny SEPT}}$ and $\alpha^{\text{\tiny IPC}}$ are nonnegative weights, with $\alpha^{\text{\tiny BR}} + \alpha^{\text{\tiny SEPT}} + \alpha^{\text{\tiny IPC}} > 0$ and $\mathcal{T}^{\text{\tiny SEPT}} \in \mathbb{T}^{\text{\tiny SEPT}}$ and $\mathcal{T}^{\text{\tiny IPC}} \in \mathbb{T}^{\text{\tiny IPC}}$.
 \end{definition}

A hybrid learning rule $\mathcal{T}$ can be rewritten as:
\begin{equation}
\label{eq:InterpHybLearnRule}
\mathcal{T} = \bar{\alpha} \left ( \frac{\alpha^{\text{\tiny BR}} }{\bar{\alpha}} \mathcal{T}^{\text{\tiny BR}} + \frac{\alpha^{\text{\tiny SEPT}} }{\bar{\alpha} } \mathcal{T}^{\text{\tiny SEPT}}+ \frac{\alpha^{\text{\tiny IPC}} }{\bar{\alpha} } \mathcal{T}^{\text{\tiny IPC}} \right ),
\end{equation} where $\bar{\alpha} := \alpha^{\text{\tiny BR}} + \alpha^{\text{\tiny SEPT}} + \alpha^{\text{\tiny IPC}}$ acts as a revision rate multiplicative factor.

There is a well-established probabilistic interpretation\footnote{See~\cite[\S4]{Sandholm2003Evolution-and-e},\cite[\S IV]{Park2019From-Population} and~\cite[\S2]{Kara_2023aa} as well as references therein.} for how a learning rule models the behavior of a strategic agent. Namely, \cite{Sandholm2010Pairwise-compar} uses (\ref{eq:InterpHybLearnRule}) to interpret $\mathcal{T}$ as modeling an agent that at every strategy revision opportunity selects a BR, SEPT or IPC learning rule with probabilities $\frac{\alpha^{\text{\tiny BR}} }{\bar{\alpha}}$, $\frac{\alpha^{\text{\tiny SEPT}} }{\bar{\alpha}}$ and $\frac{\alpha^{\text{\tiny IPC}} }{\bar{\alpha}}$, respectively.

The Lyapunov method proposed in~\cite{Hofbauer2009Stable-games-an}, or its generalization via $\delta$-passivity discussed in \S\ref{sec:DeltaPassivityAndNES}, has been studied for $\mathcal{T}^{\text{\tiny BR}}$, $\mathbb{T}^{\text{\tiny IPC}}$ and $\mathbb{T}^{\text{\tiny SEPT}}$ in isolation. The absence of a generalization of these important studies to the case of hybrid learning rules motivates the following problem, which is the central focus of this paper.

\begin{itemize}
\item[] {\bf Main Problem:} {\it Characterize convex cones of hybrid learning rules for which $\delta$-passivity holds. }
\end{itemize}

When stating our Main Problem, we recognize that it would be too \underline{ambitious} or perhaps even \underline{impossible} to establish that the entire convex cone in Definition~\ref{def:HybLearnRule} is $\delta$-passive. Notably, stability and passivity may not be preserved under linear combinations with positive coefficients. This is the case even for linear systems as the convex combination of two Hurwitz stable matrices is not, in general, Hurwitz stable. Consequently, in this article, our goal is to characterize as many and as large convex cones of $\delta$-passive hybrid learning rules as we can, which we will do in \S\ref{sec:MainResults}.

Solving the Main Problem brings the following benefits:
\begin{itemize}
    \item[(i)] The larger the convex cone of $\delta$-passive learning rules the more general the results in \S\ref{subsubsec:deltapassivitycontractivegames}-\ref{subsec:DeltaPassForDesign} will be. Furthermore, the existence of a storage function for a learning rule guarantees the existence of a global Lyapunov function as discussed in \S\ref{subsubsec:deltapassivitycontractivegames} and \S\ref{subsec:DeltaPassivityUsefulPDM}. Such a Lyapunov function can also be used to show the existence of the invariant sets required in~\cite{Benaim2003Deterministic-a} to characterize the stationary behavior of the population state when the number of agents is finite but large.
    \item[(ii)] The design methods outlined in \S\ref{subsec:DeltaPassForDesign} typically do not require knowledge of the learning rule so long as it is $\delta$-passive. Hence, these methods are all the more general the larger the convex cone of $\delta$-passive hybrid learning rules is. Furthermore, the existence of a storage function for the learning rule may be used to show the existence of a global Lyapunov function that can be used to establish global asymptotic convergence to desired equilibria and also derive anytime bounds as is done in~\cite{Martins2023Epidemic-popula}.
\end{itemize}

\section{Main Results}
\label{sec:MainResults}

In this section, we state three theorems that partially solve the {\it main problem} formulated in \S\ref{sec:HybridLearnRuleMainProb}.

\begin{theorem} \label{thm:BR-IPC} Consider that $\mathcal{T}$ is a hybrid learning rule expressible as 
\begin{equation}
\label{eq:ConvexCone1}
\mathcal{T} = \alpha^{\text{\tiny BR}} \mathcal{T}^{\text{\tiny BR}} + \alpha^{\text{\tiny IPC}}\mathcal{T}^{\text{\tiny IPC}},
\end{equation}
with $\alpha^{\text{\tiny BR}} + \alpha^{\text{\tiny IPC}}>0$, where $\mathcal{T}^{\text{\tiny IPC}}$ is in $\mathbb{T}^{\text{\tiny IPC}}$. For any number $n \geq 2$ of strategies, the hybrid learning rule $\mathcal{T}$ is $\delta$-passive with storage function $${\mathcal{S} = \alpha^{\text{\tiny BR}} \mathcal{S}^{\text{\tiny BR}} + \alpha^{\text{\tiny IPC}}\mathcal{S}^{\text{\tiny IPC}}},$$ and $\mathcal{P}$ as required in Definition~\ref{def:deltapassivity} is obtained as $${\mathcal{P} = (\alpha^{\text{\tiny BR}})^2 \mathcal{P}^{\text{\tiny BR}} + (\alpha^{\text{\tiny IPC}})^2\mathcal{P}^{\text{\tiny IPC}}}.$$
\end{theorem}

Theorem~1 ascertains that the entire convex cone formed by the combination of the best response learning rule with any IPC learning rule is $\delta$-passive. However, it is important to note that this does not include a SEPT learning rule and thus does not cover the full range of hybrid rules defined in Definition~\ref{def:HybLearnRule}.

\begin{example} We can use Theorem~1 to conclude that the following learning rules chosen arbitrarily in the convex cone characterized by~(\ref{eq:ConvexCone1}) are $\delta$-passive:
\begin{itemize}
\item[] $\mathcal{T}_{ij}(x,p)=\mathcal{T}^{\text{\tiny BR}}_{ij}(x,p)+e^{j^2[p_j-p_i]_+^3}-1+3j[p_j-p_i]_+^4$,
\item[]
\item[] $\mathcal{T}_{ij}(x,p)=2\mathcal{T}^{\text{\tiny BR}}_{ij}(x,p)+j[p_j-p_i]_++3[p_j-p_i]_+^2$~.
\end{itemize}
\end{example}

The following theorem permits the incorporation of SEPT learning rules, although it restricts IPC learning rules to those of the ND-IPC type.

{
\begin{theorem} \label{thm:BR-SEPT-ND} Consider that $\mathcal{T}$ is a hybrid learning rule expressible as \begin{equation} \label{eq:ConvexCone2} \mathcal{T} = \alpha^{\text{\tiny BR}} \mathcal{T}^{\text{\tiny BR}} + \alpha^{\text{\tiny SEPT}}\mathcal{T}^{\text{\tiny SEPT}}+ \alpha^{\text{\tiny ND}}\mathcal{T}^{\text{\tiny ND}},
\end{equation}
with $\alpha^{\text{\tiny BR}}+\alpha^{\text{\tiny SEPT}}+\alpha^{\text{\tiny ND}}>0$, where $\mathcal{T}^{\text{\tiny SEPT}}$ is in $\mathbb{T}^{\text{\tiny SEPT}}$ and $\mathcal{T}^{\text{\tiny ND}}$ is in $\mathbb{T}^{\text{\tiny ND}}$. For any number $n \geq 2$ of strategies, the hybrid rule $\mathcal{T}$ is $\delta$-passive with storage function $${\mathcal{S} = \alpha^{\text{\tiny BR}} \mathcal{S}^{\text{\tiny BR}} + \alpha^{\text{\tiny SEPT}}\mathcal{S}^{\text{\tiny SEPT}}}+ \alpha^{\text{\tiny ND}}\mathcal{S}^{\text{\tiny ND}},$$ and $\mathcal{P}$ as required in Definition~\ref{def:deltapassivity} is obtained as $$\mathcal{P} = (\alpha^{\text{\tiny BR}})^2 \mathcal{P}^{\text{\tiny BR}} + (\alpha^{\text{\tiny SEPT}})^2\mathcal{P}^{\text{\tiny SEPT}}+ (\alpha^{\text{\tiny ND}})^2\mathcal{P}^{\text{\tiny ND}}.$$
\end{theorem}}

\begin{example} We can use Theorem~2 to conclude that the following learning rules chosen arbitrarily in the convex cone characterized by~(\ref{eq:ConvexCone2}) are $\delta$-passive:{
\begin{itemize}
\item[] $\mathcal{T}_{ij}(x,p)=\mathcal{T}^{\text{\tiny BR}}_{ij}(x,p)+e^{j^3[\hat{p}_j]_+^2}-1+3[p_j-p_i]_+^4$,
\item[]
\item[] $\mathcal{T}_{ij}(x,p)=2\mathcal{T}^{\text{\tiny BR}}_{ij}(x,p)+j[\hat{p}_j]_+^3+4[p_j-p_i]_+^2$.
\end{itemize}}
\end{example}

For the case when $n=2$, Theorem~\ref{thm:BR-SEPT-ND} can be fully generalized as shown below to the entire convex cone of hybrid learning rules of Definition~\ref{def:HybLearnRule}.

\begin{theorem} \label{thm:BR-SEPT-IPC} Consider that $\mathcal{T}$ is a hybrid learning rule as specified in Definition~\ref{def:HybLearnRule}. For  the case where $n = 2$, the hybrid learning rule $\mathcal{T}$ is $\delta$-passive with storage function $${\mathcal{S} = \alpha^{\text{\tiny BR}} \mathcal{S}^{\text{\tiny BR}} + \alpha^{\text{\tiny SEPT}}\mathcal{S}^{\text{\tiny SEPT}}}+ \alpha^{\text{\tiny IPC}}\mathcal{S}^{\text{\tiny IPC}},$$ and $\mathcal{P}$ as required in Definition~\ref{def:deltapassivity} is obtained as $${\mathcal{P} = (\alpha^{\text{\tiny BR}})^2 \mathcal{P}^{\text{\tiny BR}} + (\alpha^{\text{\tiny SEPT}})^2\mathcal{P}^{\text{\tiny SEPT}}}+ (\alpha^{\text{\tiny IPC}})^2\mathcal{P}^{\text{\tiny IPC}}.$$
\end{theorem}

\begin{example} We can use Theorem~3 to conclude that the following learning rules chosen arbitrarily in the convex cone specified by Definition~\ref{def:HybLearnRule} are $\delta$-passive when $n=2$:
\begin{itemize}
\item[] $\mathcal{T}_{ij}(x,p)=\mathcal{T}^{\text{\tiny BR}}_{ij}(x,p)+e^{j^3[\hat{p}_j]_+^2}-1+3[p_j-p_i]_+^4$,
\item[]
{\item[] $\mathcal{T}_{ij}(x,p)=2\mathcal{T}^{\text{\tiny BR}}_{ij}(x,p)+j[\hat{p}_j]_+^3+[p_j-p_i]_+ e^{-(p_j-p_i)}$.}
\end{itemize}
\end{example}

\section{Auxiliary Results and Proofs}
\label{sec:AuxResults}

We now present some auxiliary results that will be helpful in proving Theorems \ref{thm:BR-IPC}, \ref{thm:BR-SEPT-ND}, and \ref{thm:BR-SEPT-IPC}. In particular, Lemma~\ref{lemma:BRplusA} might prove valuable in future work.

\begin{definition}[Positive Correlation]
A learning rule $\mathcal{T}$ is positive correlated (PC) if the EDM associated to it, $\mathcal{V}$, is such that for any $x \in \mathbb{X}$ and $p \in \mathbb{R}^n$ 
    \begin{align*}
        \mathcal{V}(x,p) \neq 0 \Rightarrow p'\mathcal{V}(x,p)>0.
    \end{align*}
\end{definition}

\begin{lemma} \label{lemma:BRplusA}
    Consider a hybrid learning rule
    \begin{align*}
        \mathcal{T} := \alpha \mathcal{T}^{\text{\tiny BR}} + \beta \mathcal{T}^{\text{\tiny A}},
    \end{align*}
    where $\mathcal{T}^{\text{\tiny A}}$ is learning rule that is PC and $\delta$-passive, having a differentiable storage function $\mathcal{S}^{\text{\tiny A}}$. 
    If for all $x \in \mathbb{X}$ and $p \in \mathbb{R}^n$ we have $\frac{\partial \mathcal{S}^{\text{\tiny A}}}{\partial x}(x,p) \mathcal{V}^{\text{\tiny BR}}(x,p)\leq 0$, 
    then for any $\alpha\geq0$ and $\beta\geq0$, such that $\alpha+\beta>0$, the hybrid learning rule $ \mathcal{T}$ is $\delta$-passive with
    \begin{align}
        \mathcal{S}(x,p) &= \alpha \mathcal{S}^{\text{\tiny BR}}(x,p) + \beta \mathcal{S}^{\text{\tiny A}}(x,p),\\
        \mathcal{P}(x,p) &= \alpha^2 \mathcal{P}^{\text{\tiny BR}}(x,p) + \beta^2 \mathcal{P}^{\text{\tiny A}}(x,p).
    \end{align}
\end{lemma}
\begin{proof}
    The dynamics $\mathcal{V}$ associated to $\mathcal{T}$ is 
    \begin{align*}
        \mathcal{V} = \alpha\mathcal{V}^{\text{\tiny BR}}+\beta\mathcal{V}^{\text{\tiny A}}.
    \end{align*}
     Let $s\st = \mathcal{S}(x\st,p\st)$. By Proposition \ref{prop:BRDeltaPassive} and definition \ref{def:deltapassivity}, for almost all $t$, the time-derivative of $s\st$ satisfies
        \begin{align*}
        \dot{s}\st =&  \alpha \dot{p}\st' \mathcal{V}^{\text{\tiny BR}}(x\st,p\st) - \alpha p\st'\dot{x}\st,\\ 
        &+\beta  \frac{\partial \mathcal{S}^{\text{\tiny A}}}{\partial x}(x\st,p\st) \dot{x}\st +\beta  \frac{\partial \mathcal{S}^{\text{\tiny A}}}{\partial p}(x\st,p\st) \dot{p}\st,\\
        \dot{s}\st \leq&  \dot{p}\st' \dot{x}\st - \alpha^2 \mathcal{P}^{\text{\tiny BR}}(x\st,p\st) - \beta^2 \mathcal{P}^{\text{\tiny A}}(x\st,p\st) \\ &- \alpha \beta p\st'\mathcal{V}^{\text{\tiny A}}(x\st,p\st) \\ & +\alpha \beta \frac{\partial \mathcal{S}^{\text{\tiny A}}}{\partial x}(x\st,p\st) \mathcal{V}^{\text{\tiny BR}}(x\st,p\st).
    \end{align*}
    Since $\frac{\partial \mathcal{S}^{\text{\tiny A}}}{\partial x}(x,p) \mathcal{V}^{\text{\tiny BR}}(x,p)\leq 0$ and $\mathcal{V}^{\text{\tiny A}}$ is PC, we obtain 
    \begin{align*}
        \dot{s}\st \leq&  \dot{p}\st' \dot{x}\st - \alpha^2 \mathcal{P}^{\text{\tiny BR}}(x\st,p\st) - \beta^2 \mathcal{P}^{\text{\tiny A}}(x\st,p\st).
    \end{align*}
\end{proof}

\subsection{Proof of Theorem \ref{thm:BR-IPC}}

The storage functions in Propositions \ref{prop:IPCDeltaPassive} and \ref{prop:SEPTDeltaPassive} will satisfy the condition in Lemma \ref{lemma:BRplusA} on being decreasing along best-response directions, that is,

\begin{proposition} \label{prop:IPC_Sdecreasing_alongBR}
    For any $\mathcal{T}^{\text{\tiny IPC}}$, $\forall x \in \mathbb{X}$, $\forall p\in \mathbb{R}^n$
    $$ \frac{\partial \mathcal{S}^{\text{\tiny IPC}}}{\partial x}(x,p) \mathcal{V}^{\text{\tiny BR}}(x,p)\leq 0.$$
\end{proposition}
\begin{proposition} \label{prop:SEPT_Sdecreasing_alongBR}
    For any $\mathcal{T}^{\text{\tiny SEPT}}$, $\forall x \in \mathbb{X}$, $\forall p\in \mathbb{R}^n$
    $$ \frac{\partial \mathcal{S}^{\text{\tiny SEPT}}}{\partial x}(x,p) \mathcal{V}^{\text{\tiny BR}}(x,p)\leq 0.$$
\end{proposition}

Therefore, from Propositions \ref{prop:IPCDeltaPassive} and \ref{prop:IPC_Sdecreasing_alongBR}, and Lemma \ref{lemma:BRplusA} we obtain Theorem \ref{thm:BR-IPC}.

\subsection{Proof of Theorem \ref{thm:BR-SEPT-ND}}
{
\begin{proposition} \label{prop:SEPT_ND_DeltaPassive}
    The learning rule 
    \begin{align} \label{eq:T_SEPT_and_ND}
        \mathcal{T}^{\text{\tiny A}} = \alpha^{\text{\tiny SEPT}}\mathcal{T}^{\text{\tiny SEPT}} +\alpha^{\text{\tiny ND}}\mathcal{T}^{\text{\tiny ND}}
    \end{align}
    is $\delta$-passive with
\begin{align}
    \mathcal{S}^{\text{\tiny A}} =& \alpha^{\text{\tiny SEPT}}\mathcal{S}^{\text{\tiny SEPT}} +\alpha^{\text{\tiny ND}}\mathcal{S}^{\text{\tiny ND}},\label{eq:S_SEPT_ND} \\ 
    \mathcal{P}^{\text{\tiny A}} =& (\alpha^{\text{\tiny SEPT}})^2\mathcal{P}^{\text{\tiny SEPT}}+ (\alpha^{\text{\tiny ND}})^2\mathcal{P}^{\text{\tiny ND}}.
\end{align}
\end{proposition}
\begin{proof}
    See the \nameref{appendix:T_SEPT_and_ND}.
\end{proof}}

Since \eqref{eq:S_SEPT_ND} is a conic combination of $\mathcal{S}^{\text{\tiny ND}}$ and $\mathcal{S}^{\text{\tiny SEPT}}$, from Propositions \ref{prop:IPC_Sdecreasing_alongBR} and \ref{prop:SEPT_Sdecreasing_alongBR} we obtain that
\begin{align}
    \frac{\partial \mathcal{S}^{\text{\tiny A}}}{\partial x}(x,p) \mathcal{V}^{\text{\tiny BR}}(x,p)&\leq 0, &\forall x \in \mathbb{X}, \forall p\in \mathbb{R}^n. \label{eq:A_Sdecreasing_alongBR}
\end{align}

Proposition \ref{prop:SEPT_ND_DeltaPassive}, \eqref{eq:A_Sdecreasing_alongBR}, and Lemma \ref{lemma:BRplusA} imply Theorem \ref{thm:BR-SEPT-ND}.

\subsection{Proof of Theorem \ref{thm:BR-SEPT-IPC}}

\begin{proposition} \label{prop:n:2_SEPT_IPC_DeltaPassive}
    For two strategies, $n=2$, the learning rule $$\mathcal{T}^{\text{\tiny A}} = \alpha^{\text{\tiny SEPT}}\mathcal{T}^{\text{\tiny SEPT}} +\alpha^{\text{\tiny IPC}}\mathcal{T}^{\text{\tiny IPC}}$$ is $\delta$-passive with
\begin{align}
    \mathcal{S}^{\text{\tiny A}} =& \alpha^{\text{\tiny SEPT}}\mathcal{S}^{\text{\tiny SEPT}} +\alpha^{\text{\tiny IPC}}\mathcal{S}^{\text{\tiny IPC}},\label{eq:n2_S_SEPT_IPC} \\ 
    \mathcal{P}^{\text{\tiny A}} =& (\alpha^{\text{\tiny SEPT}})^2\mathcal{P}^{\text{\tiny SEPT}}+ (\alpha^{\text{\tiny IPC}})^2\mathcal{P}^{\text{\tiny IPC}}.
\end{align}
\end{proposition}
\begin{proof}
    The derivative of the storage function along solutions of \eqref{eq:EDM-DEF} will be 
    \begin{align*}
        \dot{\mathcal{S}}^{\text{\tiny A}}(x\st,p\st) &= 
        \alpha^{\text{\tiny SEPT}}\dot{\mathcal{S}}^{\text{\tiny SEPT}}(x\st,p\st) \\ &+\alpha^{\text{\tiny IPC}}\dot{\mathcal{S}}^{\text{\tiny IPC}}(x\st,p\st), \\
        \dot{\mathcal{S}}^{\text{\tiny A}}(x\st,p\st) &\leq \dot{x}\st' \dot{p}\st + \alpha^{\text{\tiny SEPT}} \alpha^{\text{\tiny IPC}} \mathcal{H}(x\st,p\st) \\
        &- (\alpha^{\text{\tiny SEPT}})^2\mathcal{P}^{\text{\tiny SEPT}}(x\st,p\st) \\ &- (\alpha^{\text{\tiny IPC}})^2\mathcal{P}^{\text{\tiny IPC}}(x\st,p\st),
    \end{align*}
    where $\mathcal{H}(x,p)$ is $$\frac{\partial \mathcal{S}^{\text{\tiny SEPT}}}{\partial x}(x,p) \mathcal{V}^{\text{\tiny IPC}}(x,p) + \frac{\partial \mathcal{S}^{\text{\tiny IPC}}}{\partial x}(x,p) \mathcal{V}^{\text{\tiny SEPT}}(x,p).$$
    Moreover, for any $x$ in $\mathbb{X}$ and $p$ in $\mathbb{R}^n$, from Proposition~\ref{prop:SEPTDeltaPassive} and any IPC rule, as they are PC, we obtain  
    \begin{align} 
        \frac{\partial \mathcal{S}^{\text{\tiny SEPT}}}{\partial x}(x,p) \mathcal{V}^{\text{\tiny IPC}}(x,p) =& \nonumber \\
        \left(-\sum_{j=1}^n \phi_j([\hat{p}_j]_+) \right)& p' \mathcal{V}^{\text{\tiny IPC}}(x,p) \leq 0, \label{eq:S_SEPT_V_IPC}
    \end{align}
    which holds for any number of strategies $n$. For $n=2$,
    \begin{align*}
    \frac{\partial \mathcal{S}^{\text{\tiny IPC}}}{\partial x}(x,p) \mathcal{V}^{\text{\tiny SEPT}}(x,p) &= \int_0^{\tilde{p}_{12}} \phi_2(\tau) d\tau \mathcal{V}_1^{\text{\tiny SEPT}}(x,p) \\ &+ \int_0^{-\tilde{p}_{12}} \phi_1(\tau) d\tau \mathcal{V}_2^{\text{\tiny SEPT}}(x,p). 
    \end{align*}
    Without loss of generality assume that $\tilde{p}_{12}>0$, then, due to SEPT being PC we obtain that $\mathcal{V}^{\text{\tiny SEPT}}(x,p)\neq 0$ implies $\mathcal{V}_1^{\text{\tiny SEPT}}(x,p) < 0$, leading to
    \begin{equation*}
    \resizebox{1.0\linewidth}{!}{$
    \begin{aligned}
    \frac{\partial \mathcal{S}^{\text{\tiny IPC}}}{\partial x}(x,p) \mathcal{V}^{\text{\tiny SEPT}}(x,p) &= \underbrace{\int_0^{\tilde{p}_{12}} \phi_2(\tau) d\tau}_{\geq0} \mathcal{V}_1^{\text{\tiny SEPT}}(x,p),\leq 0.
     \end{aligned}
    $}
    \end{equation*}
    Such that $\mathcal{H}(x,p)\leq 0$ for all $x$ in $\mathbb{X}$ and $p$ in $\mathbb{R}^n$. Thus, we obtain
    \begin{align*}
        \dot{\mathcal{S}}^{\text{\tiny A}}(x\st,p\st) &\leq \dot{x}\st' \dot{p}\st - (\alpha^{\text{\tiny SEPT}})^2\mathcal{P}^{\text{\tiny SEPT}}(x\st,p\st) \\ &\quad - (\alpha^{\text{\tiny IPC}})^2\mathcal{P}^{\text{\tiny IPC}}(x\st,p\st),
    \end{align*}
    finalizing the proof.
\end{proof}

The proof of Theorem \ref{thm:BR-SEPT-IPC} follows from Propositions \ref{prop:IPC_Sdecreasing_alongBR}, \ref{prop:SEPT_Sdecreasing_alongBR}, and \ref{prop:n:2_SEPT_IPC_DeltaPassive}, and Lemma \ref{lemma:BRplusA}.

\begin{figure}[t]
    \centering
     \scalebox{0.80}{
    \begin{tikzpicture}[outer sep=auto]
    \node [anchor=east] (source) at (-3.7,0.0)  {S};
    \node [] (mode1a) at (-1.6,2) {};
    \node [] (mode2a) at (-1.6,0) {};
    \node [] (mode3a) at (-1.6,-2) {};
    \node [] (mode1b) at (1.6,2) {};
    \node [] (mode2b) at (1.6,0) {};
    \node [] (mode3b) at (1.6,-2) {};
    \node [anchor=west] (target) at (3.3,0.0)  {D};
    \draw[line width=1pt, ->] (source) |- (mode1a);
    \draw[line width=1pt, ->] (source) -- (mode2a);
    \draw[line width=1pt, ->] (source) |- (mode3a);
    \draw[line width=1pt, ->] (mode1a) -- (mode1b);
    \draw[line width=1pt, ->] (mode2a) -- (mode2b);
    \draw[line width=1pt, ->] (mode3a) -- (mode3b);
    \draw[line width=1pt, ->] (mode1b) -| (target);
    \draw[line width=1pt, ->] (mode2b) -- (target);
    \draw[line width=1pt, ->] (mode3b) -| (target);
    \node[] () at (-1.6,3.0)  {Commute};
    \node[] () at (mode1a) {\includesvg[scale=0.04]{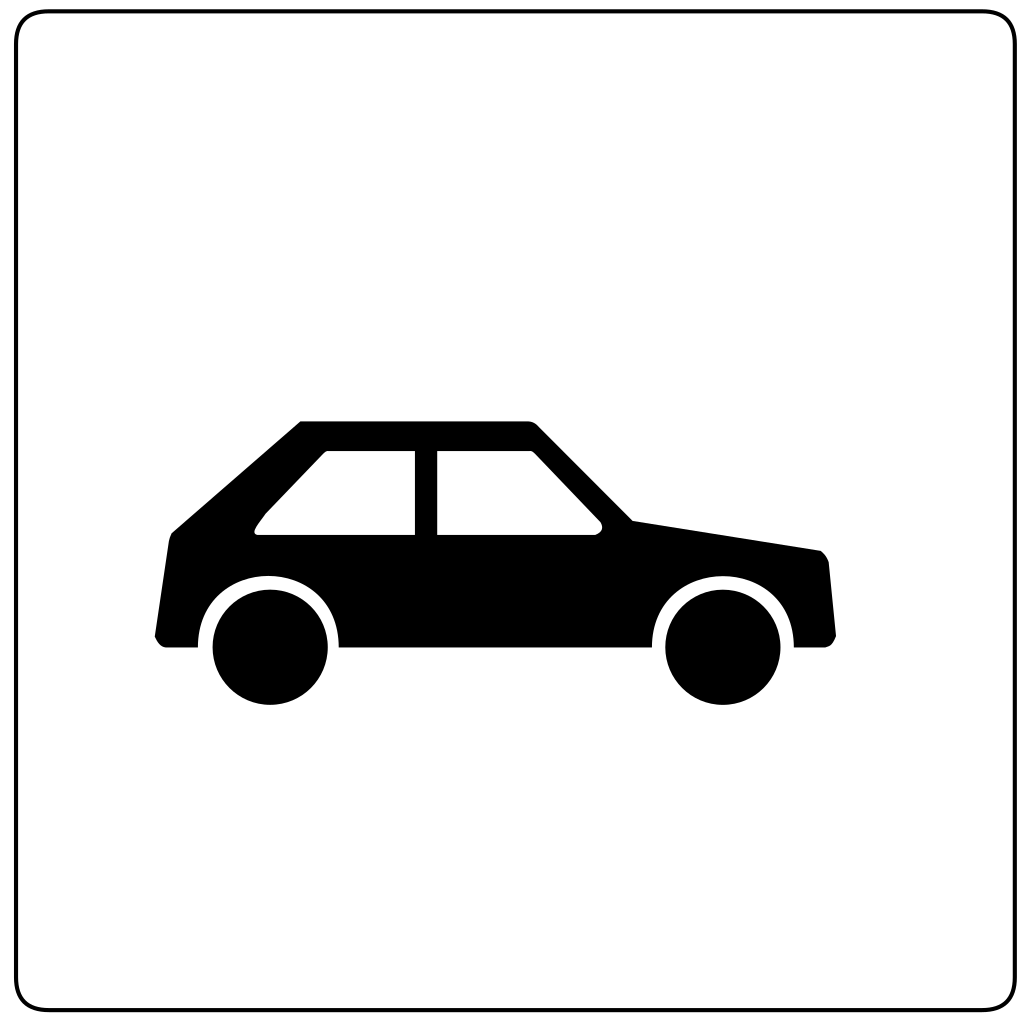} };
    \node [anchor=north, color=BrickRed] at (-1.7,1.5) {\small Main Road: $g_1(x_1+\alpha x_2)$};
    \node[] () at (mode2a) {\includesvg[scale=0.04]{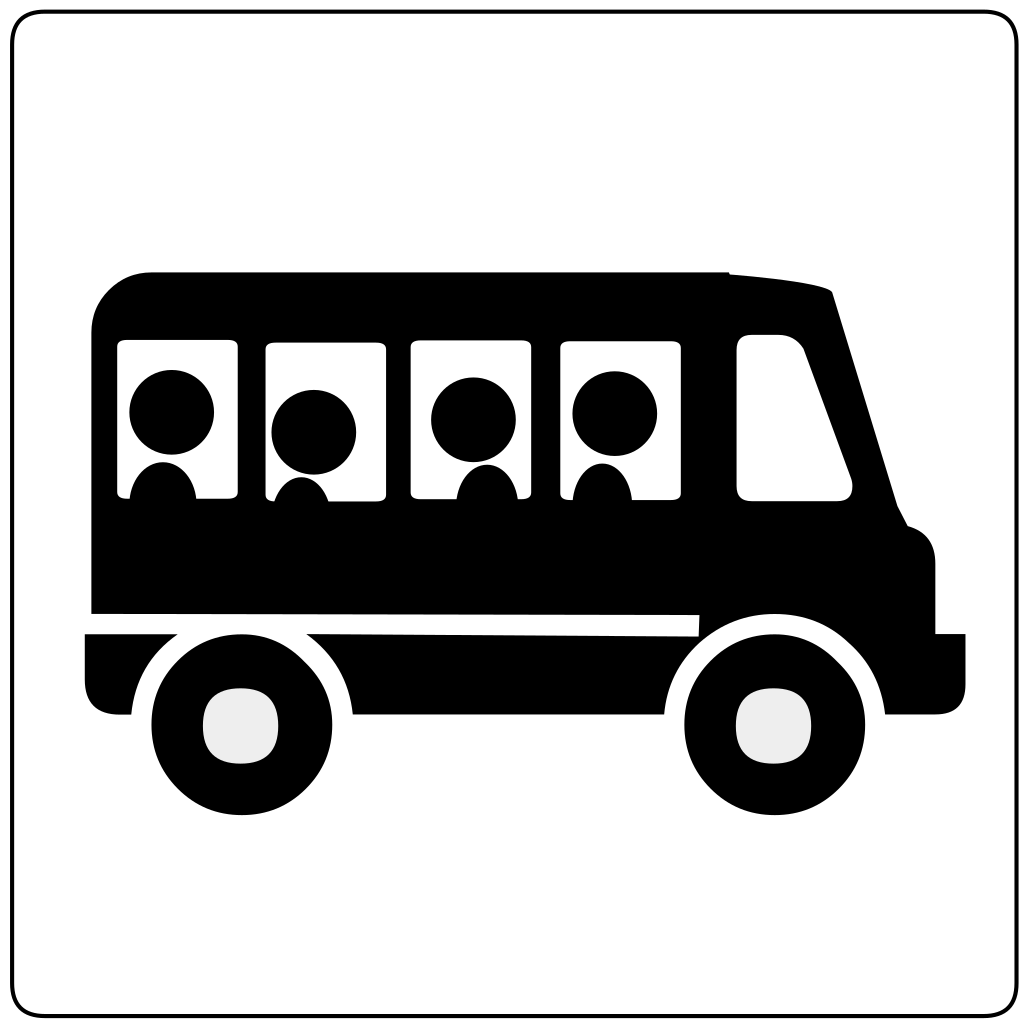} };
    \node [anchor=north, color=BrickRed] at (-1.7,-0.5) {\small Main Road: $g_1(x_1+\alpha x_2)$};
    \node[] () at (mode3a) {\includesvg[scale=0.04]{figures/car.svg} };
    \node [anchor=north, color=MidnightBlue] at (-1.6,-2.5) {\small Alternate Road: $g_2(x_3)$};
    \node[] () at (1.6,3.0)  {Mobility};
    \node[] () at (mode1b) {\includesvg[scale=0.04]{figures/car.svg} };
    \node [anchor=north, color=Violet] at (1.6,1.5) {\small Driving: $r(x_1+x_3)$};
    \node[] () at (mode2b) {\includesvg[scale=0.04]{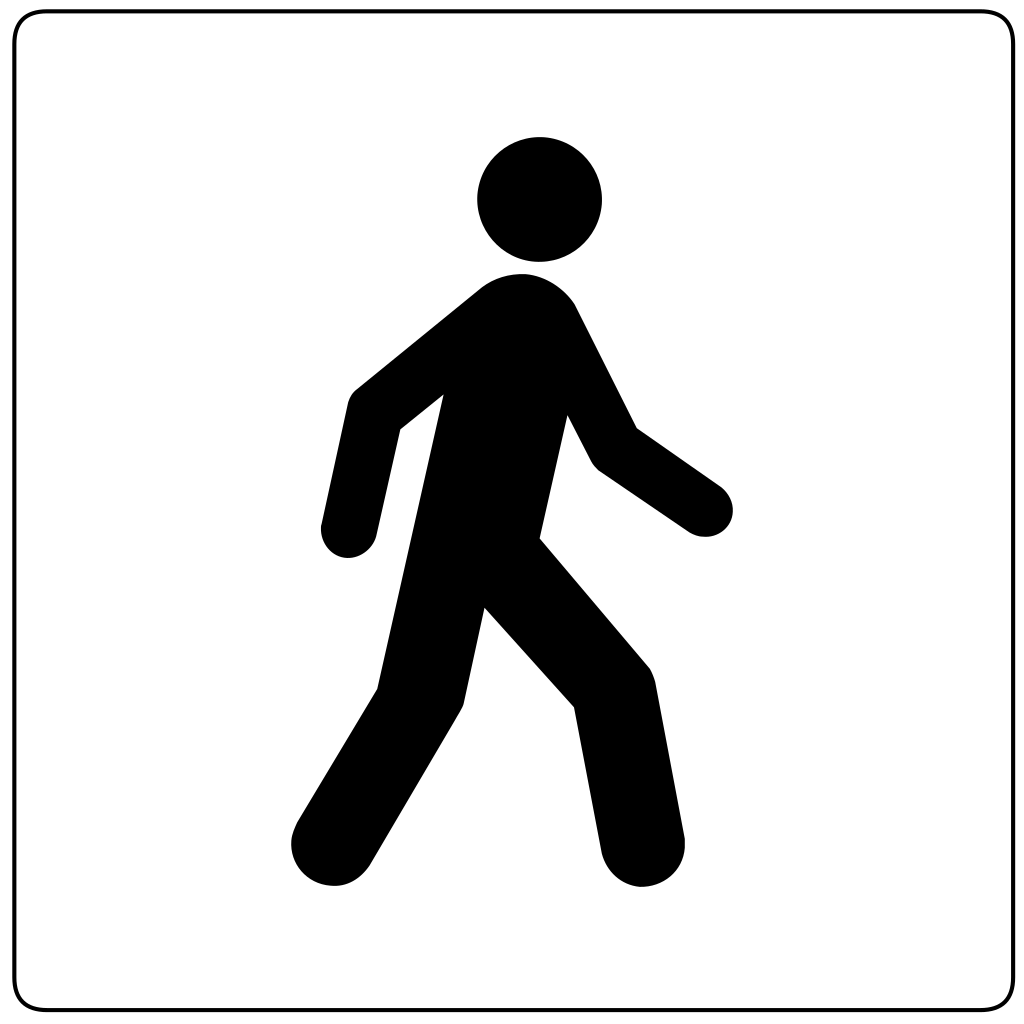} };
    \node [anchor=north, color=Black] at (1.6,-0.5) {\small Walking: $c$};
    \node[] () at (mode3b) {\includesvg[scale=0.04]{figures/car.svg} };
    \node [anchor=north, color=Violet] at (1.6,-2.5) {\small Driving: $r(x_1+x_3)$};
    \end{tikzpicture}
     }
    \caption{Available commuting strategies, resources used by each of them, and associated latency functions. `S' and `D' denote, respectively, the suburb and downtown.}
    \label{fig:strategies_paths}
\end{figure}

\begin{figure*}[t]
    \centering
    \begin{subfigure}{1.00\columnwidth}
    \resizebox{0.85\linewidth}{!}{\input{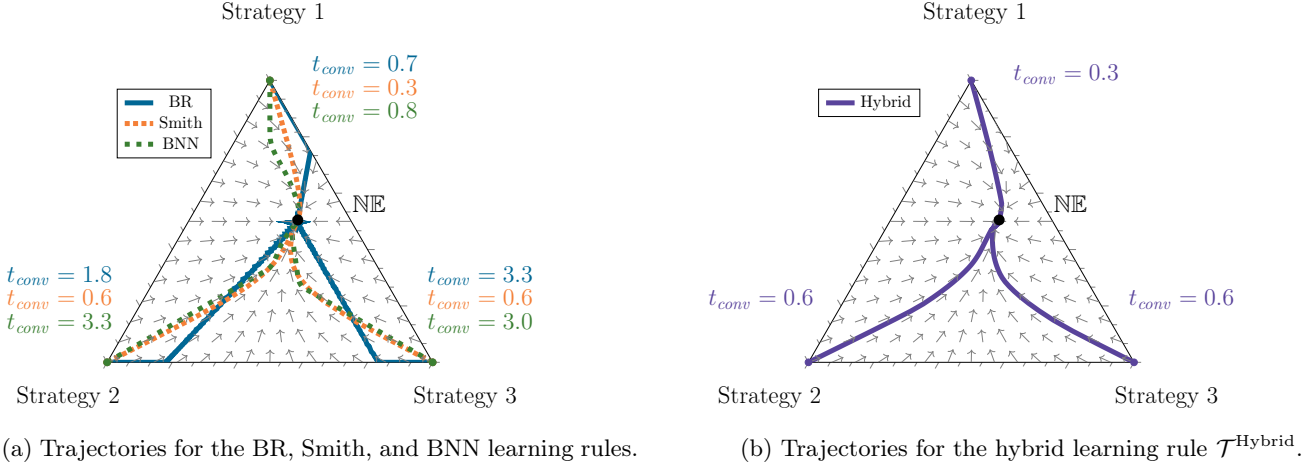}}
    \caption{Trajectories for the BR, Smith, and BNN learning rules.} \label{fig_basic}
    \end{subfigure} \hfill
    \begin{subfigure}{1.00\columnwidth}
    \resizebox{0.85\linewidth}{!}{\input{simulation/hybrid_rule.tex}}
    \caption{Trajectories for the hybrid learning rule $\mathcal{T}^{\text{\tiny Hybrid}}$.} \label{fig_hybrid}
    \end{subfigure}
    \caption{Simulation results ($t_\mathit{conv}$ indicates the earliest time at which the condition $\|x(t) - x(T)\|_1 < 10^{-3}$ holds).} \label{fig_result}
\end{figure*}

\section{Simulations} \label{sec:sims}

We illustrate our findings via simulations on a weighted congestion game, described in Example~\ref{example_game}, exploring various initial conditions and employing different learning rules, including hybrid rules.

\begin{example} \label{example_game}

As shown in Fig.~\ref{fig:strategies_paths}, we consider a weighted congestion game, where the agents select their commuting mode and route for traveling from a suburb, denoted by `S', to downtown, denoted by `D'. The strategies include: 1) driving through the main road and then on the streets downtown; 2) taking the bus along the main road and then walking in downtown; 
3) driving through an alternate road, and then driving again on the streets downtown. Each strategy utilizes a set of resources (e.g., car or bus) and undergoes a certain amount of latency. The game map is defined as follows:
\begin{align}
    \mathcal{F}(x) := -\begin{bmatrix}
        g_1(x_1+\alpha x_2) + r(x_1+x_3) \\
        g_1(x_1+\alpha x_2) + c \\
        g_2(x_3) + r(x_1+x_3)
    \end{bmatrix},
\end{align}
where ${g_1, g_2, r: \mathbb{R}_{\geq0} \rightarrow \mathbb{R}}$ are nondecreasing differentiable functions, ${\alpha \in (0,1)}$, and $c>0$. Given a population state $x$, $g_1(x_1+\alpha x_2)$, $g_2(x_3)$, $r(x_1+x_3)$, and $c$ are, respectively, the latency associated with driving on the main road, driving on the alternate road, driving on the streets downtown, and walking around downtown. Buses induce congestion at a reduced rate $\alpha$, as they can accommodate more passengers per unit area than cars.
\end{example}
The game models commuting from a suburb to the downtown area, encompassing the journey to downtown as well as moving within downtown. The agents select their transportation mean, car or bus. If they use their cars, they can also choose whether to use the main road or an alternate road to downtown, where the alternate road is likely to be faster, being used by fewer agents, but prone to higher congestion.
The payoff for each strategy is determined by the total latency experienced by a user following it. We assume the simulation horizon to be small, such that the number of buses or the width of the roads remain constant.

The game in Example~\ref{example_game} does not fall within the category of potential games, for which the population state will globally asymptotically converge to the set ${\mathbb{NE}(\mathcal{F})}$ under any learning rule that satisfies both PC and Nash stationarity conditions~\cite{Sandholm2001Potential-games}.
The game cannot be potential because the Jacobian matrix of the payoff function ${\mathcal{F}(x)}$, denoted as $D\mathcal{F}(x)$, is asymmetric~\cite{Kara2023Excess-Payoff-E}.
However, if $\cF$ is contractive, satisfying
\begin{align*}
      & \cA(x) := W(\cD\cF(x)'+\cD\cF(x))W' \preceq 0, \\
      & W := \begin{bmatrix}
    1 & -1 & 0 \\ 0 & 1 & -1
\end{bmatrix},
\end{align*}
then a population following a learning rule that satisfies Definition~\ref{def:deltapassivity} will converge to ${\mathbb{NE}(\mathcal{F})}$, as discussed in \S\ref{subsubsec:deltapassivitycontractivegames}.

Since the latency functions are nondecreasing, $\cA(x)\preceq 0$ holds if the following inequality holds:
\begin{align}
\begin{split} \label{eq:to_be_contractive}
     4 \frac{dr}{dz}(x_1+x_3) &\left(\frac{dg_2}{dz}(x_3) +\frac{dg_1}{dz} (x_1+\alpha x_2)\right) \\ &\geq \left((1-\alpha)\frac{dg_1}{dz}(x_1+\alpha x_2)\right)^2.
\end{split}
\end{align}
In particular, \eqref{eq:to_be_contractive} will be satisfied if
\begin{align}
    \label{eq:suff_cond_to_be_contractive}
    4\frac{dr}{dz}(v) \geq \frac{dg_1}{dz}(u) \text{ holds for any } v,u \in [0,1], 
\end{align}
which is the case whenever the main road has a larger traffic flow capacity than the streets downtown.  

In this example, we adopt an affine model for the latency functions, motivated by the related
literature~\cite{brown2016avoiding}, and carry out the simulations over a time horizon ${T := 10}$. Consistently with the constraints posed by the problem, we select ${c = 35}$, $g_1(z) = 10 z + 20$, $g_2(z) = 30 z + 15$, $r(z) = 30 z + 10$, and $\alpha = 1/20$. The resulting game is contractive, since the latency functions satisfy \eqref{eq:suff_cond_to_be_contractive}. Our implementation is available online~\cite{certorio_2024_11508806}.

Figure~\ref{fig_result} represents the trajectories obtained for different learning rules and initial conditions, along with the convergence time for each trajectory.  The gray arrows show the directions of the projections of the payoff vectors $\mathcal{F}(x)$ onto $\mathbb{X}$ at various points, computed as ${(I - \frac{1}{3} \mathbf{1}\mathbf{1}')\mathcal{F}(x)}$, for $x$ in the interior of $\mathbb{X}$.
For any learning rule $\cT$ as in Definition~\ref{def:HybLearnRule}, with EDM $\cV$, the direction of the arrows in Fig.~\ref{fig_result} correlates with the direction of the trajectories, since $p'\cV(x,p)\geq 0$ holds for any $x \in \mathbb{X}$ and $p \in \mathbb{R}^n$.  
The convergence time $t_{\mathit{conv}}$ is defined by the earliest time $t$ such that ${\|x(t) - x(T)\|_1 < 10^{-3}}$.

As shown in Fig.~\ref{fig_basic}, the BR, Smith, and BNN learning rules effectively guide the population towards Nash equilibria, given their $\delta$-passive nature. As expected, this is also the case for 
the hybrid learning rule ${\mathcal{T}^{\text{\tiny Hybrid}} := \frac{1}{3} \left(\mathcal{T}^{\text{\tiny BR}} + \mathcal{T}^{\text{\tiny Smith}} + \mathcal{T}^{\text{\tiny BNN}}\right)}$, generating the trajectory in Fig.~\ref{fig_hybrid}, since it also satisfies the $\delta$-passivity condition. 

\section{Conclusions}
\label{sec:conclusion}

We investigated if $\delta$-passivity is preserved under conic combinations of $\delta$-passive learning rules, and determined convex cones of hybrid rules where all elements satisfy $\delta$-passivity. Our findings can be summarized  as follows:
i) we extended the notion of $\delta$-passivity to apply it to discontinuous learning rules;
ii) we established that the best-response learning rule is $\delta$-passive;
iii) we demonstrated that hybrid rules combining best-response with impartial pairwise comparison are $\delta$-passive; and 
iv) we showed that hybrid rules combining best-response, separable excess payoff target, and {nondecreasing impartial pairwise comparison learning rules are $\delta$-passive.}

Theorems~\ref{thm:BR-IPC}, \ref{thm:BR-SEPT-ND}, and \ref{thm:BR-SEPT-IPC} also impact prior research on Nash equilibrium seeking, the analysis of payoff dynamic models, and the use of $\delta$-passivity in designing dynamic payoff mechanisms.
Our results ensure that any previous work employing $\delta$-passive rules in the population context applies to both the best-response rule and the hybrid rules we identified as $\delta$-passive.
Our simulations demonstrate that a hybrid learning rule, which combines BR, ND-IPC, and BNN rules, effectively guides the population to Nash equilibria in a congestion game. 

Future work will investigate whether any combination of BR, SEPT, and IPC can be proven to be $\delta$-passive, or if there exist hybrid rules that are not.
{Specifically, it remains to be determined if all rules satisfying  Definition~\ref{def:HybLearnRule}, are $\delta$-passive with a storage function of the form
\begin{equation*} \mathcal{S} = \alpha^{\text{\tiny BR}} \mathcal{S}^{\text{\tiny BR}} + \alpha^{\text{\tiny SEPT}}\mathcal{S}^{\text{\tiny SEPT}}+ \alpha^{\text{\tiny IPC}}\mathcal{S}^{\text{\tiny IPC}}.
\end{equation*}}

\appendix
\section*{Appendix}\label{appendix:T_SEPT_and_ND}
\stepcounter{section}
The dynamics associated to \eqref{eq:T_SEPT_and_ND} is
\begin{align*}
        \mathcal{V}^{\text{\tiny A}} = \alpha^{\text{\tiny SEPT}} \mathcal{V}^{\text{\tiny SEPT}}+ \alpha^{\text{\tiny ND}} \mathcal{V}^{\text{\tiny ND}},
\end{align*}
and the following holds for the derivative of the storage function along the trajectories:
\begin{multline}
    \dot{\mathcal{S}}^{\text{\tiny A}}(x\st,p\st) \leq
    -\mathcal{P}^{\text{\tiny A}}(x\st,p\st) + \dot{x}\st' \dot{p}\st \\
    + \alpha^{\text{\tiny SEPT}}\alpha^{\text{\tiny ND}}\frac{\partial \mathcal{S}^{\text{\tiny SEPT}}}{\partial x}(x\st,p\st) \mathcal{V}^{\text{\tiny ND}}(x\st,p\st) \\
    + \alpha^{\text{\tiny SEPT}}\alpha^{\text{\tiny ND}}\frac{\partial \mathcal{S}^{\text{\tiny ND}}}{\partial x}(x\st,p\st) \mathcal{V}^{\text{\tiny SEPT}}(x\st,p\st).   \label{eq:intermediate_1}
\end{multline}
Let $J(x,p) := 2\frac{\partial \mathcal{S}^{\text{\tiny ND}}}{\partial x}(x,p) \mathcal{V}^{\text{\tiny SEPT}}(x,p)$. We have
\begin{align*}
     J(x,p) = h(p)' \left( \phi([\hat{p}]_+)-\mathbf{1}' \phi([\hat{p}]_+) x\right) ,    
\end{align*}
where $\phi$ is as in Definition \ref{def:SEPT},  $[\hat{p}]_+$ is the vector whose $i$-th element is $[\hat{p}_i]_+$, and $h:~\mathbb{R}^n~\rightarrow~\mathbb{R}_{\geq0}^n$, for each $k$, {is defined by 
    $h_k(p) := \lambda  \sum_{j=1}^{n} \int_{0}^{[\tilde{p}_{ij}]_+}  \psi_j(s) ds,$
where $\psi_j$ is as described in Definition~\ref{def:ND}.}

We consider some arbitrary $p$ and relabel the strategy indices such that $p_i \leq p_{i+1}$ for $i=1,\dots,n-1$. This relabeling affects the order of the elements of $p,x,h,$ and $\phi$, but can be done without loss of generality. {
We define the following indices and sets to aid in the main argument:
$k(x) := \min\{ j \in [n] \ | \ p_j \geq p' x \},$
$$\overline{K}(x) := \{ j \in [n] \ | \ p_j > p_{k(x)} \}, \quad \overline{k}(x) := \min \overline{K}(x),$$
$$\underline{K}(x) := \{ j \in [n] \ | \ p_j < p_{k(x)} \}, \quad \underline{k}(x) := \max \underline{K}(x).$$}
With the vector $p$ ordered, $h_k(p)$ is nonincreasing in $k$, and $h_k(p)= h_{k+1}(p)$ if $p_k=p_{k+1}$. That, along with $\phi_j(0)=0$ as in Definition \ref{def:SEPT}, leads to
{\begin{align*}
    J(x,p) & \leq h_{k(x)}(p) \mathbf{1}' \phi([\hat{p}]_+) - \mathbf{1}' \phi([\hat{p}]_+) h(p)' x,\\
    & \leq \underbrace{\mathbf{1}' \phi([\hat{p}]_+)}_{\geq 0 } (h_{k(x)}(p) - h(p)' x).
\end{align*}
We then formulate the following problems
\begin{align}
    \text{min}_{z} \{ h(p)' z \ | \ z \in \mathbb{X}, k(z) = k(x)\},\label{op:grad_average} \\
    \text{min}_{z} \{ h(p)' z \ | \ z \in \mathbb{X}, p_{k(x)} \geq p' z\}.\label{op:grad_average_relax}
\end{align}
Where \eqref{op:grad_average_relax} relaxes \eqref{op:grad_average}, and has a Lagrangian
\begin{equation*}
\resizebox{1.0\linewidth}{!}{$
\begin{aligned}
& \mathcal{L}(z,\gamma,\mu,\theta) := h(p)' z - \gamma' z + \theta (p' z-p_{k(x)}) +\mu(\mathbf{1}' z-1).
 \end{aligned}
$}
\end{equation*}

The following proposition is used to complete the proof, as we also show that its condition holds.}
\begin{proposition}
    \begin{align}
    \max_{i \in \overline{K}(x)} \frac{h_{k(x)}(p) -h_i(p)}{p_i - p_{k(x)}} \leq \min_{i \in \underline{K}(x)} \frac{h_{k(x)}(p) -h_i(p)}{p_i - p_{k(x)}} \label{eq:minmax}\end{align}
    implies that $\mu = -h_{k(x)}(p) - \theta p_{k(x)}$, 
    \begin{align*}
    \theta &= \max_{i \in \overline{K}(x)} \frac{h_{k(x)}(p) -h_i(p)}{p_i - p_{k(x)}}, \\
    \gamma_i &= \begin{cases}
         0   & \text{ for } p_i=p_{k(x)}\\
         h_i(p)+\theta p_i + \mu   & \text{ otherwise}
    \end{cases} &\text{for each } i,
\end{align*}
    and $\Tilde{z}$, with $\Tilde{z}_{k(x)}=1$ and $\Tilde{z}_{i}=0$ if $i \neq k(x)$,
    satisfy the KKT conditions. Therefore, $\Tilde{z}$ is a solution of \eqref{op:grad_average_relax}, which is a relaxation of \eqref{op:grad_average}, and $J(x,p) \leq 0$.
\end{proposition}

To show that \eqref{eq:minmax} holds, we first note that, if either $\overline{K}(x)$ or $\underline{K}(x)$ are empty, the claim holds trivially, so we focus on the case $\overline{K}(x)\neq \emptyset$ and $\underline{K}(x)\neq \emptyset$. 
{To prove \eqref{eq:minmax} we start by introducing the proposition:
\begin{proposition} \label{pro:1}
For any $j \in [n]$, the following expression is non-decreasing with respect to $x_1$ and $x_2$:
    $$\int_{x_1}^{x_2} \frac{ \psi_j(s)}{x_2 - x_1} ds, \quad \text{with $x_1 \leq x_2$.}$$
\end{proposition}

\textbf{Step 1}:
To show that inequality \eqref{eq:minmax} holds, we begin by establishing the intermediate bound:
\begin{align*}
  \min_{i \in \underline{K}(x)} \frac{h_{k(x)}(p) -h_i(p)}{p_i - p_{k(x)}} = \frac{h_{k(x)}(p) -h_{\underline{k}(x)}(p)}{p_{\underline{k}(x)} - p_{k(x)}}.
 \end{align*}
We compute expressions for $\frac{h_{k(x)}(p) - h_i(p)}{p_i - p_{k(x)}}$ for $i \in \underline{K}(x)$:
\begin{equation*}
\resizebox{1.0\linewidth}{!}{$
\begin{aligned}
& \frac{h_{k(x)}(p) - h_i(p)}{p_i - p_{k(x)}} = 2 \lambda \sum_{j=1}^{n} \int_{{[\tilde{p}_{i j}]_+}}^{[\tilde{p}_{k(x)j}]_+} \frac{\psi_j(s)}{p_{k(x)} - p_i} ds = 2 \lambda \left( \sum_{j=k(x)}^{n} \int_{\tilde{p}_{ij}}^{\tilde{p}_{k(x)j}} \frac{\psi_j(s)}{p_i - p_{k(x)}} ds  \right). \\
\end{aligned}
$}
\end{equation*}
Applying Proposition~\ref{pro:1}, we can bound as follows:
\begin{equation*}
\resizebox{1.0\linewidth}{!}{$
\begin{aligned}
& \min_{i \in \underline{K}(x)} \frac{h_{k(x)}(p) -h_i(p)}{p_i - p_{k(x)}} = \min_{i \in \underline{K}(x)} 2 \lambda \left( \sum_{j=i}^{n} \int_{\tilde{p}_{k(x)j}}^{\tilde{p}_{ij}} \frac{\psi_j(s)}{p_{k(x)} - p_i} ds \right) \\
& = 2 \lambda \left( \sum_{j={\underline{k}(x)}}^{n} \int_{\tilde{p}_{k(x)j}}^{\tilde{p}_{{\underline{k}(x)}j}} \frac{\psi(s)}{p_{k(x)} - p_{\underline{k}(x)}} ds \right) = \frac{h_{k(x)}(p) - h_{\underline{k}(x)}(p)}{p_{\underline{k}(x)} - p_{k(x)}}.
 \end{aligned}
$}
\end{equation*}

\textbf{Step 2}:
Following Step 1, we establish the intermediate bound for $i \in \overline{K}(x)$:
\begin{align*}
 \max_{i \in \overline{K}(x)} \frac{h_{k(x)}(p) -h_i(p)}{p_i - p_{k(x)}} = \frac{h_{k(x)}(p) -h_{\overline{k}(x)}(p)}{p_{\overline{k}(x)} - p_{k(x)}}.
 \end{align*}
We compute expressions for $\frac{h_{k(x)}(p) - h_i(p)}{p_i - p_{k(x)}}$ for $i \in \overline{K}(x)$:
\begin{equation*}
\resizebox{1.0\linewidth}{!}{$
\begin{aligned}
& \frac{h_{k(x)}(p) - h_i(p)}{p_i - p_{k(x)}} = 2 \lambda \sum_{j=1}^{n} \int_{{[\tilde{p}_{i j}]_+}}^{[\tilde{p}_{k(x)j}]_+} \frac{\psi_j(s)}{p_i - p_{k(x)}} ds = 2 \lambda \left( \sum_{j=k(x)}^{n} \int_{\tilde{p}_{ij}}^{\tilde{p}_{k(x)j}} \frac{\psi_j(s)}{p_i - p_{k(x)}} ds  \right). \\
\end{aligned}
$}
\end{equation*}
Applying Proposition~\ref{pro:1}, we can bound as follows:
\begin{equation*}
\resizebox{1.0\linewidth}{!}{$
\begin{aligned}
& \max_{i \in \overline{K}(x)} \frac{h_{k(x)}(p) -h_i(p)}{p_i - p_{k(x)}} = \max_{i \in \overline{K}(x)} 2 \lambda \left( \sum_{j=k(x)}^{n} \int_{\tilde{p}_{ij}}^{\tilde{p}_{k(x)j}} \frac{\psi_j(s)}{p_i - p_{k(x)}} ds  \right) \\
& = 2 \lambda \left( \sum_{j=k(x)}^{n} \int_{\tilde{p}_{\overline{k}(x)j}}^{\tilde{p}_{k(x)j}} \frac{\psi_j(s)}{p_{\overline{k}(x)} - p_{k(x)}} ds  \right) = \frac{h_{k(x)}(p) -h_{\overline{k}(x)}(p)}{p_{\overline{k}(x)} - p_{k(x)}}.
 \end{aligned}
$}
 \end{equation*}

To complete the proof, we now show that
\begin{equation*}
\resizebox{1.0\linewidth}{!}{$
\begin{aligned}
  & \max_{i \in \overline{K}(x)} \frac{h_{k(x)}(p) -h_i(p)}{p_i - p_{k(x)}} = \frac{h_{k(x)}(p) -h_{\overline{k}(x)}(p)}{p_{\overline{k}(x)} - p_{k(x)}} \\
  & = 2 \lambda \left( \sum_{j=k(x)}^{n} \int_{\tilde{p}_{\overline{k}(x)j}}^{\tilde{p}_{k(x)j}} \frac{\psi_j(s)}{p_{\overline{k}(x)} - p_{k(x)}} ds  \right) \leq 2 \lambda \left( \sum_{j={\underline{k}(x)}}^{n} \int_{\tilde{p}_{k(x)j}}^{\tilde{p}_{{\underline{k}(x)}j}} \frac{\psi(s)}{p_{k(x)} - p_{\underline{k}(x)}} ds \right) \\
  & = \frac{h_{k(x)}(p) -h_{\underline{k}(x)}(p)}{p_{\underline{k}(x)} - p_{k(x)}} = \min_{i \in \underline{K}(x)} \frac{h_{k(x)}(p) -h_i(p)}{p_i - p_{k(x)}},
 \end{aligned}
$}
\end{equation*}
by the inequalities $(\tilde{p}_{\overline{k}(x)j} \leq \tilde{p}_{k(x)j} \leq \tilde{p}_{{\underline{k}(x)}j})$ and Proposition~\ref{pro:1}. 
Thus completing the proof of Proposition \ref{prop:SEPT_ND_DeltaPassive}.}

\bibliographystyle{elsarticle-num}        %
\bibliography{MartinsRefs,LocalRefs}

\end{document}

%% file: simulation/hybrid_rule.tex
\begin{tikzpicture}
\begin{ternaryaxis}[
every axis x label/.style={at={(ticklabel cs:1.1)},anchor=near ticklabel, xshift=-1.0cm},
every axis y label/.style={at={(ticklabel cs:1.1)},anchor=near ticklabel, xshift=1.0cm, yshift=-1.0cm},
every axis z label/.style={at={(ticklabel cs:1.1)},anchor=near ticklabel},
xlabel=\large{Strategy 1},
ylabel=\large{Strategy 2},
zlabel=\large{Strategy 3},
xticklabel=\empty,
yticklabel=\empty,
zticklabel=\empty,
minor tick num=2,
axis line style={draw=none},
grid=none,
legend pos=north west]
\addplot3 [mark = none, line width=3.0pt, color = Violet, solid] coordinates{(1.0, 0.0, 0.0)
(0.9987276815207601, 0.00033191574311652763, 0.0009404027361233825)
(0.9862438248951725, 0.0035894787045900357, 0.010166696400237506)
(0.9193404345676283, 0.021102879047664905, 0.05955668638470688)
(0.863938255665481, 0.0357552509086201, 0.10030649342589898)
(0.8173738803225249, 0.0482690627563038, 0.13435705692117134)
(0.7777465581328875, 0.059142528667968125, 0.16311091319914442)
(0.7436634216206464, 0.0687289570346345, 0.18760762134471926)
(0.7140804469907363, 0.07728496800119614, 0.20863458500806772)
(0.6881996549383353, 0.08500087076302125, 0.22679947429864347)
(0.6654178366370986, 0.09199695845249359, 0.2425852049104078)
(0.6453837098952774, 0.0981890499770011, 0.25642724012772156)
(0.6276722335462237, 0.10369756314809225, 0.2686302033056841)
(0.6119085021070567, 0.10866591404589844, 0.27942558384704486)
(0.5977981285768363, 0.11319729234841634, 0.28900457907474736)
(0.5851061246519623, 0.11736690889777138, 0.29752696645026633)
(0.5736481899282467, 0.12123158643738434, 0.3051202236343689)
(0.563494407820737, 0.12487575815576653, 0.3116298340234965)
(0.5546087862147995, 0.12979876460752338, 0.31559244917767715)
(0.5469090390976318, 0.1337389705348646, 0.3193519903675036)
(0.5400471552341887, 0.13796321009968224, 0.32198963466612907)
(0.5339873260286687, 0.1417867370575993, 0.3242259369137319)
(0.5286191516771891, 0.14524844570103232, 0.3261324026217785)
(0.5238506616348606, 0.1483831514293783, 0.32776618693576104)
(0.5196046060108405, 0.1512221493558609, 0.3291732446332986)
(0.5158156218318224, 0.15379364582883176, 0.33039073233934585)
(0.5124280391904653, 0.15612310568026008, 0.33144885512927463)
(0.509394162542053, 0.15823354149236502, 0.33237229596558193)
(0.5066729096877006, 0.16014576124083998, 0.3331813290714595)
(0.504228723390526, 0.16187858466090865, 0.3338926919485654)
(0.5020411501936644, 0.1634317314575481, 0.3345271183487877)
(0.49992152475074697, 0.16514755316259638, 0.33493092208665676)
(0.49803199520620195, 0.1665321524650383, 0.33543585232875983)
};\addlegendentry{Hybrid}
\addplot3 [color = Violet, mark = *, nodes near coords=\large{$t_{\mathit{conv}} = 0.3$}, every node near coord/.append style={xshift=2.0cm, yshift=-0.2cm}, forget plot] coordinates {(1.0, 0.0, 0.0)};
\addplot3 [mark = none, line width=3.0pt, color = Violet, solid, forget plot] coordinates{(0.0, 1.0, 0.0)
(0.0006230814503665752, 0.9985980495565843, 0.0007788689930491914)
(0.006747018815603549, 0.9848172938212131, 0.008435687363183406)
(0.05892992989861083, 0.8673451736078764, 0.07372489649351277)
(0.11196570058441474, 0.7483275631900634, 0.1397067362255219)
(0.1509763546193466, 0.6617021061412017, 0.18732153923945177)
(0.18114484849401574, 0.595858378852389, 0.2229967726535953)
(0.20691965980919574, 0.544109966508818, 0.24897037368198632)
(0.23026950137142826, 0.5022115398019014, 0.2675189588266705)
(0.25048157602928955, 0.46754384373927915, 0.2819745802314314)
(0.2682800028976078, 0.4383564902922123, 0.29336350681018)
(0.28416611302666517, 0.4134217808091746, 0.30241210616416025)
(0.2984985931664348, 0.39185356238409, 0.3096478444494752)
(0.3115416083404306, 0.3729964336339532, 0.3154619580256163)
(0.32349478508908125, 0.3563553852002086, 0.3201498297107102)
(0.3345157398351544, 0.3415529701509195, 0.3239312900139261)
(0.344851939202898, 0.32840707812355086, 0.32674098267355123)
(0.3546557371288438, 0.31669817524794464, 0.3286460876232116)
(0.36396168420301434, 0.3061682988286077, 0.32987001696837803)
(0.37279925461764496, 0.2966200110742942, 0.33058073430806095)
(0.38119457316908995, 0.28789969152369366, 0.33090573530721645)
(0.38917143434168827, 0.2798860463423965, 0.3309425193159153)
(0.3967519097071579, 0.2724820290074712, 0.33076606128537095)
(0.40395671319349064, 0.26560904882197633, 0.3304342379845331)
(0.4107672563333772, 0.2592198117734399, 0.330012931893183)
(0.4170824708559004, 0.2533162899690814, 0.32960123917501827)
(0.42294136146273653, 0.24784732336744364, 0.3292113151698199)
(0.4283865185327523, 0.24276606591961364, 0.32884741554763414)
(0.4334555222025706, 0.23803255476933982, 0.3285119230280897)
(0.43818168410166375, 0.23361244453133126, 0.32820587136700513)
(0.44259465182402485, 0.22947601528932202, 0.32792933288665327)
(0.4467209069231531, 0.22559738732906326, 0.3276817057477838)
(0.4505841791977948, 0.22195389374391591, 0.3274619270582894)
(0.45420579435571695, 0.21852557483430549, 0.3272686308099777)
(0.45760496805651435, 0.21529476736102449, 0.3271002645824613)
(0.46079905635200136, 0.21224576832632502, 0.32695517532167373)
(0.4638037703363365, 0.20936455779822202, 0.3268316718654416)
(0.46663336116498366, 0.20663856887503962, 0.3267280699599769)
(0.46930078034836903, 0.20405649556298172, 0.3266427240886494)
(0.4701239448084617, 0.20160920284722322, 0.32826685234431524)
(0.47268176260763395, 0.19934030070318695, 0.32797793668917924)
(0.47359427688228, 0.19718752608273396, 0.32921819703498617)
(0.4760184980419082, 0.1951699165915007, 0.3288115853665912)
(0.4766114447362147, 0.19324861455840334, 0.3301399407053821)
(0.478928915826392, 0.19144879411475632, 0.3296222900588518)
(0.47944155881212963, 0.1897233285092624, 0.3308351126786081)
(0.4812163117921473, 0.1881077142875526, 0.33067597392030024)
(0.48183522208802554, 0.18657023068619294, 0.33159454722578163)
(0.4845911390089664, 0.1851143180538916, 0.3302945429371421)
(0.4848579244946681, 0.183728211928337, 0.33141386357699504)
(0.48522510502309507, 0.1823917545514709, 0.3323831404254341)
(0.4857456740916407, 0.18115238937184283, 0.33310193653651654)
(0.4873379921677528, 0.1799757059510909, 0.3326863018811564)
(0.4877713069055433, 0.1788450380020528, 0.333383655092404)
(0.488228368418519, 0.17778030886905455, 0.3339913227124265)
(0.49019025848849473, 0.17676586334936586, 0.3330438781621395)
(0.4904375675869315, 0.17575589657319374, 0.3338065358398748)
(0.4907901194158226, 0.17482525031588017, 0.33438463026829723)
(0.49116373898251525, 0.17394681952986543, 0.33488944148761934)
(0.492604252425747, 0.17310656849573575, 0.3342891790785173)
(0.4928828726511308, 0.17229068430855227, 0.334826443040317)
(0.49318405523881065, 0.17152026004313792, 0.3352956847180515)
(0.4956134403876823, 0.17077823811860612, 0.3336083214937116)
(0.4955588395725706, 0.17006210688969325, 0.3343790535377362)
(0.4955048666743533, 0.16925939248451402, 0.3352357408411328)
(0.4957017435391116, 0.16861284846635458, 0.33568540799453395)
(0.4959185525362669, 0.16800110494947462, 0.33608034251425867)
(0.49614456464932405, 0.16741907303850015, 0.33643636231217594)
};
\addplot3 [color = Violet, mark = *, nodes near coords=\large{$t_{\mathit{conv}} = 0.6$}, every node near coord/.append style={xshift=-1.0cm, yshift=1.0cm}, forget plot] coordinates {(0.0, 1.0, 0.0)};
\addplot3 [mark = none, line width=3.0pt, color = Violet, solid, forget plot] coordinates{(0.0, 0.0, 1.0)
(0.0006288855334219711, 0.0007739658898423502, 0.9985971485767356)
(0.006801605634177275, 0.00836543638574208, 0.9848329579800806)
(0.058839943293977154, 0.07178250931995786, 0.869377547386065)
(0.11050678787978327, 0.13273495661617146, 0.7567582555040453)
(0.14985736386292733, 0.1732855136564035, 0.6768571224806692)
(0.18243575917865293, 0.200193732615349, 0.617370508205998)
(0.2093210044898015, 0.21918965635103096, 0.5714893391591676)
(0.232192355057532, 0.23267925525563454, 0.5351283896868335)
(0.25209398060033716, 0.242221276942949, 0.5056847424567139)
(0.26970898053682746, 0.24886934485246417, 0.4814216746107085)
(0.28556197076588113, 0.253209722580396, 0.46122830665372305)
(0.30031450597403386, 0.2550007049366323, 0.444684789089334)
(0.31419222409247427, 0.254821716080603, 0.4309860598269229)
(0.32719976973349657, 0.2534487119311892, 0.4193515183353144)
(0.33928802770718014, 0.25151679788203485, 0.4091951744107852)
(0.3505587923686022, 0.24917261822814443, 0.40026858940325355)
(0.36109902805976507, 0.24652286233769097, 0.39237810960254416)
(0.3709816386620807, 0.24365053186238225, 0.38536782947553727)
(0.38026858791274426, 0.2406205572343931, 0.3791108548528629)
(0.3890131702492626, 0.23748394476117818, 0.3735028849895595)
(0.39726167553529024, 0.23428863018035917, 0.36844969428435087)
(0.40505434565444554, 0.23110340039288518, 0.36384225395266956)
(0.4124245905561992, 0.22794987839288636, 0.35962553105091477)
(0.41927019759566647, 0.2248905680761773, 0.35583923432815656)
(0.4255923490031696, 0.22194627122621385, 0.35246137977061687)
(0.431444464817941, 0.21911403321566994, 0.34944150196638935)
(0.43687285891398536, 0.21639060940784463, 0.3467365316781703)
(0.4419179230883804, 0.21377255862546016, 0.34430951828615974)
(0.4466150751031321, 0.21125631812957152, 0.3421286067672967)
(0.450995524805401, 0.20883826325574312, 0.34016621193885616)
(0.45508689857956447, 0.20651475441536896, 0.3383983470050668)
(0.4589137524075537, 0.20428217375464605, 0.3368040738378005)
(0.46249799655058893, 0.20213695339035187, 0.3353650500590594)
(0.4658446370718719, 0.20006933754026884, 0.3340860253878595)
(0.46894328002945257, 0.19806454820746952, 0.3329921717630781)
(0.471817748319119, 0.19612362562777985, 0.3320586260531013)
(0.4744899203821999, 0.19424722492061686, 0.3312628546971834)
(0.47697908803219846, 0.19243533129640988, 0.3305855806713918)
(0.4793023039122662, 0.1906874043506603, 0.3300102917370737)
(0.47998044469940204, 0.18901293344210224, 0.3310066218584959)
(0.4817375663442895, 0.1874346463400739, 0.3308277873156368)
(0.4823387941508616, 0.1859316774942048, 0.33172952835493374)
(0.4850774216065539, 0.18450773691433509, 0.33041484147911115)
(0.4853268234887032, 0.18315134302460906, 0.33152183348668796)
(0.48568056349694927, 0.181841835131239, 0.3324776013718119)
(0.48618408211187036, 0.18062819849174316, 0.3331877193963867)
(0.487760069601336, 0.1794755886172808, 0.3327643417813834)
(0.4881772619237658, 0.1783675197638817, 0.3334552183123527)
(0.488618733529978, 0.17732398975625185, 0.3340572767137704)
(0.49012938246385046, 0.17632949392711025, 0.33354112360903954)
(0.4904795748220377, 0.17536897682155553, 0.33415144835640703)
(0.4908529410816699, 0.17446306195044797, 0.33468399696788237)
(0.49229540922486137, 0.17359740149228398, 0.3341071892828549)
(0.4927477299687864, 0.17255047939314802, 0.3347017906380658)
(0.49304708043337525, 0.17176691829845145, 0.33518600126817355)
(0.49336454522146217, 0.17102585632480125, 0.33560959845373683)
(0.4947447979507937, 0.17031373468052854, 0.334941467368678)
(0.49496049208916015, 0.16962028736279758, 0.3354192205480425)
(0.4951981030117506, 0.16896473249015426, 0.33583716449809536)
(0.49544819219276715, 0.16834263118075454, 0.33620917662647853)
};
\addplot3 [color = Violet, mark = *, nodes near coords=\large{$t_{\mathit{conv}} = 0.6$}, every node near coord/.append style={xshift=1.0cm, yshift=1.0cm}, forget plot] coordinates {(0.0, 0.0, 1.0)};
\addplot3 [black, mark = *, nodes near coords=\large{$\mathbb{NE}$}, every node near coord/.append style={xshift=1.5cm}, line width=2.5pt, forget plot] coordinates {(0.505050505050505, 0.16161616161616155, 0.33333333333333326)};
\draw [->, thin, gray ,shorten >=3pt] (axis cs: 1.0, 0.0, 0.0)--(axis cs: 0.95, -0.012499999999999983, 0.06250000000000003);
\draw [->, thin, gray ,shorten >=3pt] (axis cs: 0.9285714285714286, 0.07142857142857151, 0.0)--(axis cs: 0.8869336169370796, 0.05056638306292057, 0.06250000000000003);
\draw [->, thin, gray ,shorten >=3pt] (axis cs: 0.9285714285714286, 0.0, 0.07142857142857151)--(axis cs: 0.873015873015873, -0.00694444444444446, 0.13392857142857145);
\draw [->, thin, gray ,shorten >=3pt] (axis cs: 0.8571428571428572, 0.1428571428571429, 0.0)--(axis cs: 0.8233726958525346, 0.11412730414746554, 0.06250000000000006);
\draw [->, thin, gray ,shorten >=3pt] (axis cs: 0.8571428571428572, 0.07142857142857151, 0.07142857142857151)--(axis cs: 0.8125476614133198, 0.053523767158108894, 0.1339285714285715);
\draw [->, thin, gray ,shorten >=3pt] (axis cs: 0.8571428571428572, 0.0, 0.1428571428571429)--(axis cs: 0.7946428571428573, 0.003125000000000128, 0.20223214285714292);
\draw [->, thin, gray ,shorten >=3pt] (axis cs: 0.7857142857142857, 0.2142857142857143, 0.0)--(axis cs: 0.759359847071988, 0.1781401529280119, 0.062499999999999924);
\draw [->, thin, gray ,shorten >=3pt] (axis cs: 0.7857142857142857, 0.1428571428571429, 0.07142857142857151)--(axis cs: 0.7512536692759296, 0.11481775929549909, 0.13392857142857162);
\draw [->, thin, gray ,shorten >=3pt] (axis cs: 0.7857142857142857, 0.07142857142857151, 0.1428571428571429)--(axis cs: 0.7358075692963751, 0.058835287846481926, 0.20535714285714277);
\draw [->, thin, gray ,shorten >=3pt] (axis cs: 0.7857142857142857, 0.0, 0.2142857142857143)--(axis cs: 0.7232142857142857, 0.019531249999999986, 0.2572544642857143);
\draw [->, thin, gray ,shorten >=3pt] (axis cs: 0.7142857142857143, 0.2857142857142858, 0.0)--(axis cs: 0.6949329224075418, 0.24256707759245844, 0.06250000000000007);
\draw [->, thin, gray ,shorten >=3pt] (axis cs: 0.7142857142857143, 0.2142857142857143, 0.07142857142857151)--(axis cs: 0.6892238330975955, 0.17684759547383308, 0.13392857142857145);
\draw [->, thin, gray ,shorten >=3pt] (axis cs: 0.7142857142857143, 0.1428571428571429, 0.1428571428571429)--(axis cs: 0.6786135444743935, 0.11602931266846361, 0.20535714285714282);
\draw [->, thin, gray ,shorten >=3pt] (axis cs: 0.7142857142857143, 0.07142857142857151, 0.2142857142857143)--(axis cs: 0.652043978748524, 0.0711703069657615, 0.2767857142857141);
\draw [->, thin, gray ,shorten >=3pt] (axis cs: 0.7142857142857143, 0.0, 0.2857142857142858)--(axis cs: 0.6517857142857145, 0.046875000000000215, 0.3013392857142858);
\draw [->, thin, gray ,shorten >=3pt] (axis cs: 0.6428571428571429, 0.3571428571428572, 0.0)--(axis cs: 0.6301256613756614, 0.30737433862433866, 0.06249999999999991);
\draw [->, thin, gray ,shorten >=3pt] (axis cs: 0.6428571428571429, 0.2857142857142858, 0.07142857142857151)--(axis cs: 0.6265354868061874, 0.2395359417652412, 0.13392857142857145);
\draw [->, thin, gray ,shorten >=3pt] (axis cs: 0.6428571428571429, 0.2142857142857143, 0.1428571428571429)--(axis cs: 0.6200152146060218, 0.17462764253683533, 0.20535714285714285);
\draw [->, thin, gray ,shorten >=3pt] (axis cs: 0.6428571428571429, 0.1428571428571429, 0.2142857142857143)--(axis cs: 0.6045048701298701, 0.11870941558441561, 0.27678571428571425);
\draw [->, thin, gray ,shorten >=3pt] (axis cs: 0.6428571428571429, 0.07142857142857151, 0.2857142857142858)--(axis cs: 0.5803571428571432, 0.10209627329192586, 0.3175465838509318);
\draw [->, thin, gray ,shorten >=3pt] (axis cs: 0.6428571428571429, 0.0, 0.3571428571428572)--(axis cs: 0.6043956043956046, 0.06250000000000014, 0.3331043956043957);
\draw [->, thin, gray ,shorten >=3pt] (axis cs: 0.5714285714285714, 0.4285714285714286, 0.0)--(axis cs: 0.5649682348901098, 0.3725317651098901, 0.06249999999999998);
\draw [->, thin, gray ,shorten >=3pt] (axis cs: 0.5714285714285714, 0.3571428571428572, 0.07142857142857151)--(axis cs: 0.5632554945054944, 0.30281593406593416, 0.13392857142857148);
\draw [->, thin, gray ,shorten >=3pt] (axis cs: 0.5714285714285714, 0.2857142857142858, 0.1428571428571429)--(axis cs: 0.5602106227106227, 0.2344322344322345, 0.20535714285714285);
\draw [->, thin, gray ,shorten >=3pt] (axis cs: 0.5714285714285714, 0.2142857142857143, 0.2142857142857143)--(axis cs: 0.5532904595404595, 0.16992382617382618, 0.27678571428571436);
\draw [->, thin, gray ,shorten >=3pt] (axis cs: 0.5714285714285714, 0.1428571428571429, 0.2857142857142858)--(axis cs: 0.5221497252747249, 0.12963598901098883, 0.34821428571428525);
\draw [->, thin, gray ,shorten >=3pt] (axis cs: 0.5714285714285714, 0.07142857142857151, 0.3571428571428572)--(axis cs: 0.539588948787062, 0.1339285714285715, 0.3264824797843666);
\draw [->, thin, gray ,shorten >=3pt] (axis cs: 0.5714285714285714, 0.0, 0.4285714285714286)--(axis cs: 0.5567226890756303, 0.06250000000000017, 0.3807773109243699);
\draw [->, thin, gray ,shorten >=3pt] (axis cs: 0.5, 0.5, 0.0)--(axis cs: 0.4994877049180328, 0.4380122950819672, 0.06249999999999998);
\draw [->, thin, gray ,shorten >=3pt] (axis cs: 0.5, 0.4285714285714286, 0.07142857142857151)--(axis cs: 0.49944196428571436, 0.36662946428571436, 0.13392857142857162);
\draw [->, thin, gray ,shorten >=3pt] (axis cs: 0.5, 0.3571428571428572, 0.1428571428571429)--(axis cs: 0.4993622448979591, 0.29528061224489804, 0.20535714285714282);
\draw [->, thin, gray ,shorten >=3pt] (axis cs: 0.5, 0.2857142857142858, 0.2142857142857143)--(axis cs: 0.4991883116883116, 0.22402597402597407, 0.2767857142857142);
\draw [->, thin, gray ,shorten >=3pt] (axis cs: 0.5, 0.2142857142857143, 0.2857142857142858)--(axis cs: 0.4985119047619048, 0.1532738095238095, 0.3482142857142858);
\draw [->, thin, gray ,shorten >=3pt] (axis cs: 0.5, 0.1428571428571429, 0.3571428571428572)--(axis cs: 0.49784482758620674, 0.20535714285714257, 0.29679802955665013);
\draw [->, thin, gray ,shorten >=3pt] (axis cs: 0.5, 0.07142857142857151, 0.4285714285714286)--(axis cs: 0.49973849372384943, 0.1339285714285716, 0.36633293484757923);
\draw [->, thin, gray ,shorten >=3pt] (axis cs: 0.5, 0.0, 0.5)--(axis cs: 0.5, 0.0625, 0.4375);
\draw [->, thin, gray ,shorten >=3pt] (axis cs: 0.4285714285714286, 0.5714285714285714, 0.0)--(axis cs: 0.43331826401446655, 0.5089285714285714, 0.05775316455696202);
\draw [->, thin, gray ,shorten >=3pt] (axis cs: 0.4285714285714286, 0.5, 0.07142857142857151)--(axis cs: 0.43451992922331906, 0.43749999999999994, 0.127980070776681);
\draw [->, thin, gray ,shorten >=3pt] (axis cs: 0.4285714285714286, 0.4285714285714286, 0.1428571428571429)--(axis cs: 0.4364639200390053, 0.3660714285714285, 0.1974646513895661);
\draw [->, thin, gray ,shorten >=3pt] (axis cs: 0.4285714285714286, 0.3571428571428572, 0.2142857142857143)--(axis cs: 0.44014550264550273, 0.2946428571428572, 0.2652116402116402);
\draw [->, thin, gray ,shorten >=3pt] (axis cs: 0.4285714285714286, 0.2857142857142858, 0.2857142857142858)--(axis cs: 0.44977678571428575, 0.2232142857142858, 0.32700892857142866);
\draw [->, thin, gray ,shorten >=3pt] (axis cs: 0.4285714285714286, 0.2142857142857143, 0.3571428571428572)--(axis cs: 0.4910714285714288, 0.17938311688311714, 0.32954545454545464);
\draw [->, thin, gray ,shorten >=3pt] (axis cs: 0.4285714285714286, 0.1428571428571429, 0.4285714285714286)--(axis cs: 0.45114087301587313, 0.1827876984126984, 0.3660714285714286);
\draw [->, thin, gray ,shorten >=3pt] (axis cs: 0.4285714285714286, 0.07142857142857151, 0.5)--(axis cs: 0.44097720746590097, 0.12152279253409917, 0.4375);
\draw [->, thin, gray ,shorten >=3pt] (axis cs: 0.4285714285714286, 0.0, 0.5714285714285714)--(axis cs: 0.437192118226601, 0.05387931034482751, 0.5089285714285713);
\draw [->, thin, gray ,shorten >=3pt] (axis cs: 0.3571428571428572, 0.6428571428571429, 0.0)--(axis cs: 0.36613952063189437, 0.5803571428571429, 0.05350333651096284);
\draw [->, thin, gray ,shorten >=3pt] (axis cs: 0.3571428571428572, 0.5714285714285714, 0.07142857142857151)--(axis cs: 0.36808755760368667, 0.5089285714285714, 0.122983870967742);
\draw [->, thin, gray ,shorten >=3pt] (axis cs: 0.3571428571428572, 0.5, 0.1428571428571429)--(axis cs: 0.37106207111665634, 0.4375, 0.1914379288833438);
\draw [->, thin, gray ,shorten >=3pt] (axis cs: 0.3571428571428572, 0.4285714285714286, 0.2142857142857143)--(axis cs: 0.376164596273292, 0.36607142857142866, 0.25776397515527955);
\draw [->, thin, gray ,shorten >=3pt] (axis cs: 0.3571428571428572, 0.3571428571428572, 0.2857142857142858)--(axis cs: 0.3869505494505496, 0.2946428571428572, 0.31840659340659344);
\draw [->, thin, gray ,shorten >=3pt] (axis cs: 0.3571428571428572, 0.2857142857142858, 0.3571428571428572)--(axis cs: 0.4196428571428571, 0.2280219780219781, 0.3523351648351648);
\draw [->, thin, gray ,shorten >=3pt] (axis cs: 0.3571428571428572, 0.2142857142857143, 0.4285714285714286)--(axis cs: 0.40772275405007374, 0.2262058173784978, 0.3660714285714286);
\draw [->, thin, gray ,shorten >=3pt] (axis cs: 0.3571428571428572, 0.1428571428571429, 0.5)--(axis cs: 0.3834061550151976, 0.17909384498480252, 0.4375);
\draw [->, thin, gray ,shorten >=3pt] (axis cs: 0.3571428571428572, 0.07142857142857151, 0.5714285714285714)--(axis cs: 0.3749519969278034, 0.1161194316436253, 0.5089285714285714);
\draw [->, thin, gray ,shorten >=3pt] (axis cs: 0.3571428571428572, 0.0, 0.6428571428571429)--(axis cs: 0.37065637065637075, 0.04898648648648649, 0.5803571428571429);
\draw [->, thin, gray ,shorten >=3pt] (axis cs: 0.2857142857142858, 0.7142857142857143, 0.0)--(axis cs: 0.29821428571428577, 0.6517857142857143, 0.05);
\draw [->, thin, gray ,shorten >=3pt] (axis cs: 0.2857142857142858, 0.6428571428571429, 0.07142857142857151)--(axis cs: 0.3006136665192393, 0.5803571428571429, 0.11902919062361797);
\draw [->, thin, gray ,shorten >=3pt] (axis cs: 0.2857142857142858, 0.5714285714285714, 0.1428571428571429)--(axis cs: 0.3041153009427122, 0.5089285714285714, 0.1869561276287165);
\draw [->, thin, gray ,shorten >=3pt] (axis cs: 0.2857142857142858, 0.5, 0.2142857142857143)--(axis cs: 0.3097052247587669, 0.4375, 0.2527947752412333);
\draw [->, thin, gray ,shorten >=3pt] (axis cs: 0.2857142857142858, 0.4285714285714286, 0.2857142857142858)--(axis cs: 0.32004527162977875, 0.3660714285714285, 0.3138832997987928);
\draw [->, thin, gray ,shorten >=3pt] (axis cs: 0.2857142857142858, 0.3571428571428572, 0.3571428571428572)--(axis cs: 0.34566326530612257, 0.29464285714285715, 0.35969387755102034);
\draw [->, thin, gray ,shorten >=3pt] (axis cs: 0.2857142857142858, 0.2857142857142858, 0.4285714285714286)--(axis cs: 0.3482142857142858, 0.268765133171913, 0.38302058111380144);
\draw [->, thin, gray ,shorten >=3pt] (axis cs: 0.2857142857142858, 0.2142857142857143, 0.5)--(axis cs: 0.32755750605326883, 0.2349424939467312, 0.43749999999999994);
\draw [->, thin, gray ,shorten >=3pt] (axis cs: 0.2857142857142858, 0.1428571428571429, 0.5714285714285714)--(axis cs: 0.3134661513859276, 0.17760527718550112, 0.5089285714285714);
\draw [->, thin, gray ,shorten >=3pt] (axis cs: 0.2857142857142858, 0.07142857142857151, 0.6428571428571429)--(axis cs: 0.3065186032630323, 0.11312425387982498, 0.5803571428571429);
\draw [->, thin, gray ,shorten >=3pt] (axis cs: 0.2857142857142858, 0.0, 0.7142857142857143)--(axis cs: 0.3023809523809525, 0.04583333333333335, 0.6517857142857143);
\draw [->, thin, gray ,shorten >=3pt] (axis cs: 0.2142857142857143, 0.7857142857142857, 0.0)--(axis cs: 0.22972336416581024, 0.7232142857142857, 0.04706235011990408);
\draw [->, thin, gray ,shorten >=3pt] (axis cs: 0.2142857142857143, 0.7142857142857143, 0.07142857142857151)--(axis cs: 0.23239319092122832, 0.6517857142857143, 0.1158210947930575);
\draw [->, thin, gray ,shorten >=3pt] (axis cs: 0.2142857142857143, 0.6428571428571429, 0.1428571428571429)--(axis cs: 0.23615016872890887, 0.5803571428571429, 0.1834926884139483);
\draw [->, thin, gray ,shorten >=3pt] (axis cs: 0.2142857142857143, 0.5714285714285714, 0.2142857142857143)--(axis cs: 0.24182808716707022, 0.5089285714285714, 0.24924334140435833);
\draw [->, thin, gray ,shorten >=3pt] (axis cs: 0.2142857142857143, 0.5, 0.2857142857142858)--(axis cs: 0.2514062076443481, 0.43749999999999994, 0.31109379235565193);
\draw [->, thin, gray ,shorten >=3pt] (axis cs: 0.2142857142857143, 0.4285714285714286, 0.3571428571428572)--(axis cs: 0.27100536746490506, 0.3660714285714286, 0.36292320396366645);
\draw [->, thin, gray ,shorten >=3pt] (axis cs: 0.2142857142857143, 0.3571428571428572, 0.4285714285714286)--(axis cs: 0.27678571428571436, 0.32440476190476203, 0.39880952380952384);
\draw [->, thin, gray ,shorten >=3pt] (axis cs: 0.2142857142857143, 0.2857142857142858, 0.5)--(axis cs: 0.27377366609294324, 0.288726333907057, 0.4375000000000001);
\draw [->, thin, gray ,shorten >=3pt] (axis cs: 0.2142857142857143, 0.2142857142857143, 0.5714285714285714)--(axis cs: 0.2528314341300723, 0.23823999444135632, 0.5089285714285714);
\draw [->, thin, gray ,shorten >=3pt] (axis cs: 0.2142857142857143, 0.1428571428571429, 0.6428571428571429)--(axis cs: 0.24284174876847292, 0.17680110837438429, 0.5803571428571429);
\draw [->, thin, gray ,shorten >=3pt] (axis cs: 0.2142857142857143, 0.07142857142857151, 0.7142857142857143)--(axis cs: 0.23699357305564595, 0.11122071265863984, 0.6517857142857143);
\draw [->, thin, gray ,shorten >=3pt] (axis cs: 0.2142857142857143, 0.0, 0.7857142857142857)--(axis cs: 0.23315363881401618, 0.04363207547169811, 0.7232142857142857);
\draw [->, thin, gray ,shorten >=3pt] (axis cs: 0.1428571428571429, 0.8571428571428572, 0.0)--(axis cs: 0.1607935333896873, 0.7946428571428572, 0.04456360946745563);
\draw [->, thin, gray ,shorten >=3pt] (axis cs: 0.1428571428571429, 0.7857142857142857, 0.07142857142857151)--(axis cs: 0.16361931194339396, 0.7232142857142857, 0.11316640234232042);
\draw [->, thin, gray ,shorten >=3pt] (axis cs: 0.1428571428571429, 0.7142857142857143, 0.1428571428571429)--(axis cs: 0.167478354978355, 0.6517857142857143, 0.18073593073593078);
\draw [->, thin, gray ,shorten >=3pt] (axis cs: 0.1428571428571429, 0.6428571428571429, 0.2142857142857143)--(axis cs: 0.1730641885926188, 0.5803571428571429, 0.2465786685502384);
\draw [->, thin, gray ,shorten >=3pt] (axis cs: 0.1428571428571429, 0.5714285714285714, 0.2857142857142858)--(axis cs: 0.1818698817106461, 0.5089285714285714, 0.3092015468607826);
\draw [->, thin, gray ,shorten >=3pt] (axis cs: 0.1428571428571429, 0.5, 0.3571428571428572)--(axis cs: 0.19780680728667316, 0.4375, 0.36469319271332695);
\draw [->, thin, gray ,shorten >=3pt] (axis cs: 0.1428571428571429, 0.4285714285714286, 0.4285714285714286)--(axis cs: 0.20535714285714302, 0.38637599709934745, 0.40826686004350987);
\draw [->, thin, gray ,shorten >=3pt] (axis cs: 0.1428571428571429, 0.3571428571428572, 0.5)--(axis cs: 0.2053571428571428, 0.34369349005424954, 0.45094936708860756);
\draw [->, thin, gray ,shorten >=3pt] (axis cs: 0.1428571428571429, 0.2857142857142858, 0.5714285714285714)--(axis cs: 0.19316202090592338, 0.29790940766550533, 0.5089285714285714);
\draw [->, thin, gray ,shorten >=3pt] (axis cs: 0.1428571428571429, 0.2142857142857143, 0.6428571428571429)--(axis cs: 0.1796709410767275, 0.23997191606612975, 0.5803571428571429);
\draw [->, thin, gray ,shorten >=3pt] (axis cs: 0.1428571428571429, 0.1428571428571429, 0.7142857142857143)--(axis cs: 0.17191672229639526, 0.1762975634178906, 0.6517857142857143);
\draw [->, thin, gray ,shorten >=3pt] (axis cs: 0.1428571428571429, 0.07142857142857151, 0.7857142857142857)--(axis cs: 0.1668817093311313, 0.10990400495458306, 0.7232142857142857);
\draw [->, thin, gray ,shorten >=3pt] (axis cs: 0.1428571428571429, 0.0, 0.8571428571428572)--(axis cs: 0.16334894613583142, 0.04200819672131146, 0.7946428571428572);
\draw [->, thin, gray ,shorten >=3pt] (axis cs: 0.07142857142857151, 0.9285714285714286, 0.0)--(axis cs: 0.09151632091239806, 0.8660714285714286, 0.04241225051617343);
\draw [->, thin, gray ,shorten >=3pt] (axis cs: 0.07142857142857151, 0.8571428571428572, 0.07142857142857151)--(axis cs: 0.09442385444743943, 0.7946428571428572, 0.11093328840970357);
\draw [->, thin, gray ,shorten >=3pt] (axis cs: 0.07142857142857151, 0.7857142857142857, 0.1428571428571429)--(axis cs: 0.09829612413251286, 0.7232142857142857, 0.17848959015320157);
\draw [->, thin, gray ,shorten >=3pt] (axis cs: 0.07142857142857151, 0.7142857142857143, 0.2142857142857143)--(axis cs: 0.10370879120879127, 0.6517857142857143, 0.24450549450549453);
\draw [->, thin, gray ,shorten >=3pt] (axis cs: 0.07142857142857151, 0.6428571428571429, 0.2857142857142858)--(axis cs: 0.1118092298647855, 0.5803571428571428, 0.3078336272780718);
\draw [->, thin, gray ,shorten >=3pt] (axis cs: 0.07142857142857151, 0.5714285714285714, 0.3571428571428572)--(axis cs: 0.1252606882168927, 0.5089285714285714, 0.36581074035453603);
\draw [->, thin, gray ,shorten >=3pt] (axis cs: 0.07142857142857151, 0.5, 0.4285714285714286)--(axis cs: 0.13392857142857156, 0.4515063424947146, 0.41456508607671405);
\draw [->, thin, gray ,shorten >=3pt] (axis cs: 0.07142857142857151, 0.4285714285714286, 0.5)--(axis cs: 0.13392857142857145, 0.4040461121157324, 0.46202531645569617);
\draw [->, thin, gray ,shorten >=3pt] (axis cs: 0.07142857142857151, 0.3571428571428572, 0.5714285714285714)--(axis cs: 0.13392857142857156, 0.3564849624060152, 0.5095864661654135);
\draw [->, thin, gray ,shorten >=3pt] (axis cs: 0.07142857142857151, 0.2857142857142858, 0.6428571428571429)--(axis cs: 0.11705740578439974, 0.30258545135845766, 0.580357142857143);
\draw [->, thin, gray ,shorten >=3pt] (axis cs: 0.07142857142857151, 0.2142857142857143, 0.7142857142857143)--(axis cs: 0.10717497430626938, 0.2410393114080165, 0.6517857142857144);
\draw [->, thin, gray ,shorten >=3pt] (axis cs: 0.07142857142857151, 0.1428571428571429, 0.7857142857142857)--(axis cs: 0.1008330989876266, 0.17595261529808776, 0.7232142857142857);
\draw [->, thin, gray ,shorten >=3pt] (axis cs: 0.07142857142857151, 0.07142857142857151, 0.8571428571428572)--(axis cs: 0.09641813737181021, 0.10893900548533277, 0.7946428571428572);
\draw [->, thin, gray ,shorten >=3pt] (axis cs: 0.07142857142857151, 0.0, 0.9285714285714286)--(axis cs: 0.09316770186335412, 0.04076086956521739, 0.8660714285714286);
\draw [->, thin, gray ,shorten >=3pt] (axis cs: 0.0, 1.0, 0.0)--(axis cs: 0.02195945945945946, 0.9375, 0.04054054054054054);
\draw [->, thin, gray ,shorten >=3pt] (axis cs: 0.0, 0.9285714285714286, 0.07142857142857151)--(axis cs: 0.024899854333576128, 0.8660714285714286, 0.1090287170949954);
\draw [->, thin, gray ,shorten >=3pt] (axis cs: 0.0, 0.8571428571428572, 0.1428571428571429)--(axis cs: 0.028733221476510074, 0.7946428571428572, 0.17662392138063282);
\draw [->, thin, gray ,shorten >=3pt] (axis cs: 0.0, 0.7857142857142857, 0.2142857142857143)--(axis cs: 0.03393916913946588, 0.7232142857142857, 0.2428465451462484);
\draw [->, thin, gray ,shorten >=3pt] (axis cs: 0.0, 0.7142857142857143, 0.2857142857142858)--(axis cs: 0.04141566265060244, 0.6517857142857143, 0.3067986230636834);
\draw [->, thin, gray ,shorten >=3pt] (axis cs: 0.0, 0.6428571428571429, 0.3571428571428572)--(axis cs: 0.053062403697996904, 0.5803571428571429, 0.3665804534448603);
\draw [->, thin, gray ,shorten >=3pt] (axis cs: 0.0, 0.5714285714285714, 0.4285714285714286)--(axis cs: 0.06250000000000003, 0.5184394409937888, 0.4190605590062112);
\draw [->, thin, gray ,shorten >=3pt] (axis cs: 0.0, 0.5, 0.5)--(axis cs: 0.06250000000000003, 0.4675632911392405, 0.4699367088607595);
\draw [->, thin, gray ,shorten >=3pt] (axis cs: 0.0, 0.4285714285714286, 0.5714285714285714)--(axis cs: 0.06250000000000007, 0.41661294481691596, 0.5208870551830841);
\draw [->, thin, gray ,shorten >=3pt] (axis cs: 0.0, 0.3571428571428572, 0.6428571428571429)--(axis cs: 0.055059523809523794, 0.3645833333333334, 0.5803571428571429);
\draw [->, thin, gray ,shorten >=3pt] (axis cs: 0.0, 0.2857142857142858, 0.7142857142857143)--(axis cs: 0.04279556650246301, 0.30541871921182273, 0.6517857142857142);
\draw [->, thin, gray ,shorten >=3pt] (axis cs: 0.0, 0.2142857142857143, 0.7857142857142857)--(axis cs: 0.03502263581488937, 0.24176307847082498, 0.7232142857142858);
\draw [->, thin, gray ,shorten >=3pt] (axis cs: 0.0, 0.1428571428571429, 0.8571428571428572)--(axis cs: 0.029655612244897916, 0.17570153061224494, 0.7946428571428572);
\draw [->, thin, gray ,shorten >=3pt] (axis cs: 0.0, 0.07142857142857151, 0.9285714285714286)--(axis cs: 0.025727172312223854, 0.10820139911634766, 0.8660714285714286);
\draw [->, thin, gray ,shorten >=3pt] (axis cs: 0.0, 0.0, 1.0)--(axis cs: 0.022727272727272742, 0.039772727272727286, 0.9375);

\end{ternaryaxis}
\end{tikzpicture}